\documentclass[11pt]{article}
\usepackage{amssymb}
\usepackage{amsmath}
\usepackage{amscd}
\usepackage{amsfonts}
\usepackage{graphicx,color}
\usepackage[title]{appendix}
\usepackage{mathrsfs}

\numberwithin{equation}{section}

\newtheorem{theorem}{Theorem}[section]
\newtheorem{corollary}{Corollary}[section]
\newtheorem{lemma}{Lemma}[section]
\newtheorem{remark}{Remark}[section]
\newtheorem{proposition}{Proposition}[section]
\makeatletter
\begin{document}

\title{A novel approach to the Lindel{\"o}f hypothesis}
\author{A.S. Fokas}

\date{}
\maketitle
%\preprint{APS/123-QED}
%\title{1}
\begin{center}
Department of Applied Mathematics and Theoretical Physics,

University of Cambridge, CB3 0WA, UK,

and

Viterbi School of Engineering, University of Southern California,

Los Angeles, California, 90089-2560, USA.
\end{center}

%\noindent {\large \bf Abstract}
\begin{abstract}
Lindel{\"o}f's hypothesis,  one of the most important open problems in the history of mathematics, states that for large $t$, Riemann's zeta function $\zeta(1/2+it)$ is of
order $O(t^{\varepsilon})$ for any $\varepsilon>0$ . It is well known that for large $t$, the leading order asymptotics of the Riemann zeta function can be expressed in terms of a transcendental exponential sum. The usual approach to the Lindel\"of hypothesis involves the use of ingenious techniques for the estimation of this sum. However, since such estimates can not yield an asymptotic formula for the above sum, it appears that this strategy cannot lead to the proof of Lindel\"of's hypothesis. Here, a completely different approach is introduced. In particular, a novel linear integral equation is derived for $|\zeta(\sigma+it)|^2, \ 0<\sigma<1$ whose asymptotic analysis yields asymptotic results for a certain Riemann zeta-type double exponential sum. This sum has the same structure as the sum describing the leading asymptotics of $|\zeta(\sigma+it)|^2$, namely it involves $m_1^{-\sigma-it}m_2^{-\sigma-it}$, but its summation limits are different than those of the sum corresponding to $|\zeta(\sigma+it)|^2$. The analysis of the above integral equation requires the asymptotic estimation of four different integrals denoted by $I_1,I_2,\tilde{I}_3,\tilde{I}_4$, as well as the derivation of an exact relation between certain double exponential sums. Here, the latter relation is derived, and also the rigorous analysis of the first two integrals $I_1$ and $I_2$ is presented. For the remaining two integrals, formal results are only derived that suggest a possible roadmap for the derivation of rigorous asymptotic results of the above double exponential sum, as well as for other sums associated with $|\zeta(\sigma+it)|^2$. Additional developments  suggested by the above novel approach are also discussed.
\end{abstract}

%%%%%%%%%%%%%%%%%%%%%%%%%%%%%%%%%%%%%%
\section{Introduction}

Riemann's hypothesis, which is perhaps the most celebrated open problem in the history of   mathematics, can be verified numerically for $t$ up to
 $O(10^{13})$. This suggests that it is important to investigate the large
 $t$-asymptotics of the Riemann function $\zeta(s), \ s=\sigma+i t,
 \ 0<\sigma<1, \ t>0$. This investigation is
 closely related with Lindel{\"o}f's hypothesis. Indeed, it is well known that the leading asymptotics for large $t$ of $\zeta(s)$ can be expressed in terms of the transcendental sum
\begin{equation} \label{1.1}
\zeta(s) \sim \sum_{m=1}^{[t]} \frac{1}{m^s} , \qquad s=\sigma+it,
 \quad 0<\sigma<1,\quad t \to \infty,
\end{equation}
 where throughout this paper [A] denotes the integer part of the positive number A. Lindel{\"o}f's hypothesis, one of the most important open problems in the history of mathematics, states that for $\sigma=1/2$,
 this sum is of order $O(t^{\varepsilon})$ for any $\varepsilon>0$.

It is stated in \cite{GM} that the validity of  Lindel{\"o}f's hypothesis reduces significantly the number of possible zeros that disobey the Riemann hypothesis.
Indeed, let $N(\sigma,T)$ denote the number of zeros,
 $\beta+it$, of $\zeta(s)$, such that $\beta>\sigma$
 and $0<t\leq T$. Then,
 $$N(\sigma,T)\leq N(T)<AT \ln T,\qquad\frac{1}{2}<\sigma<1,  $$
where $N(T)$ denotes the number of zeros of $\zeta(s)$ in the domain
$0\leq \sigma\leq 1$, $t\leq T$, and $A$ is a constant. According to
Theorem 9.19A of \cite{T}, $N(\sigma,T)$ satisfies the estimate
$$N(\sigma,T)=O\left(T^{\frac{3}{2}-\sigma}(\ln T)^{5}\right),\qquad
\frac{1}{2}<\sigma<1,~~T\rightarrow\infty. $$
In \cite{HT} it is proven that if Lindel{\"o}f's hypothesis is valid then
$$N(\sigma,T)=O\left(T^{\epsilon}\right),\qquad
\frac{3}{4}+\delta\leq\sigma<1,~~T\rightarrow\infty, $$
with $\epsilon$ and $\delta$ positive and arbitrarily small. The first result in this direction was obtained in \cite{In} where under the assumption of Lindel{\"o}f's hypothesis, it was shown that
$$N(\sigma,T)=O\left(T^{2(1-\sigma)+\epsilon}\right),\qquad
\frac{1}{2}\leq\sigma\leq 1,~~T\rightarrow\infty, $$
for $\epsilon>0$, arbitrarily small.
%In \cite{HT} it is stated that the above estimate shattered for the first time ``the semblance that the Lindel{\"o}f hypothesis is much weaker than the Riemann Hypothesis and has not even an essential influence on the vertical distribution of the zeros, i.e. $N(\sigma,t)$". Furthermore, in \cite{Tu}, where  a slightly stronger form of the above estimate is proven by a different method, 
 In \cite{Tu} it is stated that Lindel{\"o}f's hypothesis  is much stronger than expected and even implies the estimate
$$N(\sigma,T)=O\left(T^{\epsilon}\right),\qquad
\frac{1}{2}+\delta\leq\sigma<1,~~T\rightarrow\infty, $$
with $\epsilon$ and $\delta$ positive and arbitrarily small.
%
%We note that the source of the above mentioned semblance, namely that the Lindel{\"o}f's hypothesis is much weaker than the Riemann Hypothesis,  is the estimate of the expression $N(T+1)-N(T)$, as $T\rightarrow\infty$, which is
%\begin{itemize}
%\item $o\left(\ln T\right)$ under the assumption of Lindel{\"o}f's hypothesis,
%\item $O\left(\dfrac{\ln T}{\ln( \ln T)}\right)$ under the assumption of the Riemann Hypothesis,
%\item $O\left(\ln T\right)$ in the unconditional case.
%\end{itemize}

%But, if Lindel{\"o}f's Hypothesis is valid, then Theorem B.5 of \cite{T}
%satisfies that
%$$N(\sigma,T+1)-N(\sigma,T)=o\left(\ln T\right),~~\frac{1}{2}<\sigma<1,~~T\rightarrow\infty,
%$$
%thus Lindel{\"o}f's Hypothesis implies a dramatic reduction of the possible
%zeros of $\zeta(s)$. See also [TH] and [BI].

 The sum of the rhs of \eqref{1.1} is a particular case of an exponential sum. Pioneering results for the estimation of such sums were obtained in 1916  using methods developed by Weyl \cite{W}, and Hardy and Littlewood \cite{HL}, when it was shown that $\zeta(1/2 + it) = O(t^{1/6 +  \varepsilon})$. In the last 80 years some slight progress was made using the ingenious techniques of Vinogradov \cite{V}. Currently, the best result is due to  Bourgain \cite{B} who has been able to reduce the exponent factor to $53/342\approx 0.155$.

It turns out that the best estimate for the growth of $\zeta(s)$ as $t \to \infty$, is based on the approximate functional equation, see page 79 of \cite{T},
\begin{multline} \label{1.2}
\zeta(s) = \sum_{n\le x} \frac{1}{n^s} + \frac{(2\pi)^s}{\pi} \sin{\left(  \frac{\pi s}{2} \right)} \Gamma (1-s) \sum_{n\le y} \frac{1}{n^{1-s}}  + O \left(  x^{-\sigma} + |t|^{\frac{1}{2}-\sigma} y^{\sigma -1}  \right),    \\
xy = \frac{t}{2\pi}, s=\sigma+i t \quad 0<\sigma<1, \quad t \to \infty,
\end{multline}
where $\Gamma(s)$, $s \in \mathbb{C}$, denotes the gamma function. It is interesting  that, in contrast to the usual situation in asymptotics  where higher order terms in an asymptotic expansion are more complicated, the higher order terms of the asymptotic expansion of  $\zeta(s)$ can be computed {\it explicitly}. Siegel, in his classical paper \cite{S}  presented the asymptotic expansion of $\zeta(s)$ to {\it all} orders in the important case of $x=y=\sqrt{t/2\pi}$. In \cite{FL}, analogous results are presented for {\it any} $x$ and $y$ valid to {\it all} orders.  A similar result for the Hurwitz zeta function is presented in \cite{FF}.

An interesting relation between the Riemann hypothesis and the solution
of a particular Neumann problem for the two-dimensional Laplace equation
is presented in \cite{FG}.

A major obstacle in trying to prove Lindel{\"o}f's hypothesis via the estimation of relevant  exponential sums is that in estimates one ``loses" something (the more powerful the technique the less the loss). Hence, it does not appear that
it is possible to prove Lindel{\"o}f's hypothesis via the above technically
complicated but conceptually straightforward approach. Actually, world renowned mathematicians believe that the only way of proving  the Lindel{\"o}f hypothesis is by
 proving the Riemann hypothesis  and  by employing the fact that   the latter 
implies the former  \cite{SA}.

\subsection*{A Novel Exact Integral Equation for $| \zeta (\sigma+it)|^2$}

Here, a different approach to Lindel{\"o}f's hypothesis is introduced. In particular, the following linear integral equation satisfied by $| \zeta (\sigma+i t) |^2$ is derived:

%the Riemann zeta function is embedded in a more complex mathematical structure and
%the asymptotics  of this new structure is computed. In this approach,  by using asymptotics, one
%avoids the loss occurring in estimates.
%
%Of course, the main difficulty with the above approach is finding a suitable
%``more complex mathematical structure".
%After many  unsuccessful attempts and with the crucial help of Jonatan Lenells and
%Anthony Ashton, it has been possible to embed the Riemann zeta
%function in the Riemann-Hilbert problem (or equivalently, singular integral equation) shown below.
%
%
%
%
% We have advocated the following approach for estimating the large
% $t$ asymptotics of the Riemann zeta and related functions: we
% first derive a suitable
% equation satisfied by the function under consideration, and
% then compute the large $t-$asymptotics of this equation. Here, we implement
% this approach for the following equation satisfied by the absolute
% value of the Riemann zeta function:

%\subsection*{Main results}

%The following equation is valid:

 \begin{multline} \label{1.3}
\frac{t}{\pi} \oint_{-\infty}^{\infty} \Re{\left\{   \frac{\Gamma(it - i\tau t)}{\Gamma(\sigma + i t)} \Gamma(\sigma + i\tau t) \right\}} \left| \zeta (\sigma+i\tau t) \right|^2 \text{d}\tau + \mathcal{G}(\sigma,t) = 0, \\
0<\sigma< 1, \quad t>0,
\end{multline}
where the principal value integral is defined with respect to $\tau=1$, and the function $\mathcal{G}(\sigma,t)$ is defined by the formulae
\begin{multline} \label{1.4}
\mathcal{G}(\sigma,t) = \begin{cases}
      \zeta(2\sigma) + \left( \frac{\Gamma(1-\bar{s})}{\Gamma(s)} + \frac{\Gamma(1-s)}{\Gamma(\bar{s})} \right) \Gamma(2\sigma-1) \zeta(2\sigma-1) + \frac{2(\sigma -1)\zeta(2\sigma-1)}{(\sigma - 1)^2 + t^2}, & \sigma \ne \frac{1}{2}, \\
     \Re{\left\{ \Psi \left(  \frac{1}{2} + it \right) \right\}} + 2\gamma - \ln{2\pi} + \frac{2}{1+4t^2}, & \sigma = \frac{1}{2},
   \end{cases}
\end{multline}
with $s=\sigma+it$, $\Psi(z)$ denoting the digamma function, i.e.
\begin{equation*}
\Psi(z) = \frac{\frac{\mathrm{d}}{\mathrm{dz}}\Gamma(z)}{\Gamma(z)}, \quad z \in \mathbb{C},
\end{equation*}
and $\gamma$ denoting the Euler constant. 

The derivation of \eqref{1.3}, which is presented in section 2,   uses the  so-called Plemelj formulae,
which are the fundamental ingredients of the theory of Riemann-Hilbert problems.

It is shown in Corollary 2.1 that $\mathcal{G}(\sigma,t)$ satisfies
\begin{equation}\label{1.4b}
\mathcal{G}(\sigma,t) = \begin{cases}
     \zeta(2\sigma) + O\left(\frac{t^{1-2\sigma}}{1-2\sigma}\right), &\sigma\neq \frac{1}{2},\\
     \ln t +O(1), &\sigma =\frac{1}{2}, \end{cases} \qquad t \to \infty.
     \end{equation}

%The  large $t$ asymptotic analysis of \eqref{1.3} yields the following asymptotic estimates:
%
%First,
%\begin{align}\label{R1}
%\hspace*{-15mm}\sum_{m_{2}=\left[t^{\frac{1}{2}-\epsilon}\right]}^{[t]}\sum_{m_{1}=1}^{\left[\frac{m_{2}}{t^{\frac{1}{2}-\epsilon}}\right]}\frac{1}{(m_{1}+m_{2})^{\frac{1}{2}+i t} m_{2}^{\frac{1}{2}-it}}  = O\left( t^{\epsilon} \ln{t} \right), \quad t \to \infty, \ \text{ for any } \epsilon>0,
%\end{align} 
% which is the analogue of the Lindel\"of hypothesis for the Riemann-type double sum defined by the lhs of the above equation.
%
%
%
%
%
%Second,
%\begin{equation}\label{R2}
%\left|\sum_{m=1}^{[t]}
%\frac{1}{m^{\frac{1}{2}+it}}\right|^2 -2\Re\left\{S^\epsilon(t)\right\}=O\left( t^{\epsilon} \ln{t} \right),  \ \quad t \to \infty, \ \text{ for any } \epsilon>0,
%\end{equation}
% with $S^\epsilon$ given by
% \begin{align}
%S^\epsilon(t)= t^{i\left(\epsilon-\frac{1}{2}\right)t}e^{i\left(t-t^{\frac{1}{2}+\epsilon}\right)}\mathop{\sum\sum}_{m_{1},m_{2}\in M^\epsilon} \frac{1}{m_{1}^{\frac{1}{2}+it}m_{2}^{\frac{1}{2}-it}}e^{-i\frac{m_{2}}{m_{1}}t^{\frac{1}{2}+\epsilon}},
%\end{align}
%and the set $M^\epsilon$ defined by
%\begin{align*}
%M^\epsilon=\left\{m_{j}=1,\ldots,[t],~j=1,2, \ \ \frac{m_{2}}{m_{1}}\in(t^{-\frac{1}{2}-\epsilon}, t^{\frac{1}{2}-\epsilon}-1)\right\}.
%\end{align*}
%Taking into consideration that $S^\epsilon(t)$ depends on $\epsilon$, \eqref{R2} suggests the validity of Lindel\"of's hypothesis.

It is expected that the asymptotic analysis of \eqref{1.3} yields useful results for the asymptotics of $| \zeta (\sigma+i t) |^2$. Regarding the analysis of \eqref{1.3} it is noted that  the term $\Gamma(it-i\tau t)\Gamma(\sigma+i\tau t)/\Gamma(\sigma+it)$, $\ -\infty<\tau<\infty$, appearing in \eqref{1.3},  decays exponentially for large $t$, unless $\tau$
is in the interval
 $$ -t^{\delta_{1}-1}\leq \tau\leq 1+t^{\delta_{4}-1}, \quad \delta_1>0, \ \delta_4>0. $$
 Thus,   equation
 \eqref{1.3} simplifies to the equation
\begin{align} \label{1.5}
\frac{t}{\pi} &\oint_{-t^{\delta_1-1}}^{1+t^{\delta_4-1}} \Re {\left\{ \frac{\Gamma(it-i\tau t)}{\Gamma(\sigma + it)}\Gamma(\sigma + i\tau t) \right\}} \left| \zeta(\sigma+i\tau t)\right|^2 \textrm{d}\tau + \mathcal{G}(\sigma,t)  \\
+ & O\left(e^{-\pi t^{\delta_{14}}}\right) = 0, \quad t\to\infty,\quad 0<\sigma<1, \  \delta_1>0, \ \delta_4>0, \ \delta_{14}=\mbox{min}(\delta_{1},\delta_{4}),\notag
\end{align}
where the principal value integral is defined with respect to
$\tau=1$.

It turns out that the computation of the large $t$ asymptotics of \eqref{1.5}
is facilitated by splitting  the interval $[-t^{\delta_{1}-1},1+t^{\delta_{4}-1}]$
into the following four subintervals:
\begin{align}\label{1.6}
&L_{1}=[-t^{\delta_{1}-1},t^{-1}],
L_{2}=[t^{-1},t^{\delta_{2}-1}], L_{3}=[t^{\delta_{2}-1},1-t^{\delta_{3}-1}],\notag\\
&L_{4}=[1-t^{\delta_{3}-1},1+t^{\delta_{4}-1}], \quad  \delta_2>0, \ \delta_3>0.
\end{align}
 Then, the asymptotic evaluation of \eqref{1.5} reduces to the analysis of the four integrals,
\begin{align}\label{1.7}
I_{j}(\sigma,t)=&\frac{t}{\pi} \oint_{L_{j}} \Re {\left\{ \frac{\Gamma(it-i\tau t)}{\Gamma(\sigma + it)}\Gamma(\sigma + i\tau t) \right\}} \left| \zeta(\sigma+i\tau t)\right|^2 \textrm{d}\tau ,\notag \\
& 0<\sigma<1, \ \  t>0, \ \ j=1,2,3,4,
\end{align}
where $I_{1}$, $I_{2}$, $I_{3}$, $I_{4}$ also depend on $\delta_{1}$,
$\delta_{2}$, $(\delta_{2},\delta_{3})$, $(\delta_{3},\delta_{4})$, respectively,   $\{L_{j}\}_{1}^{4}$
are defined in \eqref{1.6}, and the principal value integral is needed only for $I_{4}$.

It turns out that it is possible to estimate $I_1$ and $I_2$ directly. However, for $I_3$ and $I_4$ it is necessary to replace $| \zeta (\sigma+i t) |^2$ by its large $t$ asymptotics, i.e. by the sum 

\begin{equation}\label{asym_R}
S_R(\sigma,t)=\left(\sum_{m_1=1}^{[t]} \frac{1}{m_1^s}\right) \left(\sum_{m_2=1}^{[t]} \frac{1}{m_2^{\bar{s}}}\right), \qquad s=\sigma+i t, \ 0<\sigma<1, \ t>0.
\end{equation}

The  large $t$-asymptotic estimates obtained by replacing in $I_3$ and $I_4$ $|\zeta(s)|^2$   with $S_R$ can be used to obtain rigorous results via two different ways \cite{F}: one  is to use the results of \cite{FL} to estimate the error terms; in this approach the estimation of the error terms of $I_4$ is straightforward,but for $I_3$, since the error term depends on $\tau t\in\left(t^{\delta_2},t-t^{\delta_3}\right)$  $\delta_2$ cannot be arbitrarily small. An alternative approach is based on the fact that the sum $S_R$ satisfies a linear integral equation similar to \eqref{1.3}, see \cite{F}; by employing this equation there is no need to keep track of the error terms of replacing  $|\zeta(s)|^2$   with $S_R$  in $I_3$ and $I_4$, but it is necessary to handle such error terms regarding $I_1$ and $I_2$.

%Fortunately, this can be done \textit{without} the need to keep track of the relevant error terms: it is shown in \cite{F} that the sum $S_R$ satisfies a linear integral equation similar to \eqref{1.3}, thus the rigorous version of our results which will be presented in \cite{F} is based on this alternative integral equation derived in \cite{F}.

Replacing  in the definition of $I_3$ and $I_4$ $| \zeta (\sigma+i t) |^2$ by $S_R$ we obtain two integrals denoted by $\tilde{I}_3$ and $\tilde{I}_4$. The integral  $\tilde{I}_3$ is defined by
%Let $\tilde{I}_3$ denote the integral obtained from the rhs of equation \eqref{1.7} with $j=3$, when $|\zeta|^2$ is replaced by its leading term asymptotics:
\begin{align}\label{1.10}
\tilde{I}_3(\sigma,t,\delta_2, \delta_3)= \sum_{m_1=1}^{[t]}\sum_{m_2=1}^{[t]} \frac{1}{m_1^\sigma}\frac{1}{m_2^\sigma} \Re\left\{ J_3\left(\sigma,t,\delta_2, \delta_3, \frac{m_2}{m_1}\right) \right\}, 
\end{align}
where
\begin{multline} \label{1.11}
J_3\left(\sigma,t,\delta_2, \delta_3, \frac{m_2}{m_1}\right) = \frac{t}{\pi} \int_{t^{\delta_2 -1}}^{1-t^{\delta_3 -1}} \frac{\Gamma(it-it\tau)}{\Gamma(\sigma+it)} \Gamma(\sigma+it\tau) \left(\frac{m_2}{m_1}\right)^{i\tau t }  d\tau, \\ 0 < \sigma <1, \quad t>0; \quad 0 < \delta_2 <1, \  0< \delta_3 <1, \ m_j = 1,  \dots, [t], \quad j= 1,2.
\end{multline}

Similarly, after some simple calculations presented in section 6, $\tilde{I}_4$ takes the form
\begin{align}  \label{1.12a}
&\tilde{I_4}(\sigma, t,\delta_3,\delta_4) =\Re\left\{  \sum_{m_1=1}^{[t]}  \sum_{m_2=1}^{[t]} \frac{1}{m_{1}^s}\frac{1}{m_{2}^{\bar{s}}} J_4\left(\sigma, t,\delta_{2},\delta_{3},\frac{m_{1}}{m_{2}}\right)\right\},
\end{align}
where $J_4$ is defined by
\begin{align} \label{1.13a}
&J_4\left(\sigma,t,\delta_3,\delta_4,\frac{m_{1}}{m_{2}}\right) =\frac{1}{\pi}\oint_{-t^{\delta_4}}^{t^{\delta_3}} \Gamma (i x) \frac{\Gamma (\sigma + i t - i x)}{\Gamma(\sigma + i t )}\left(\frac{m_{1}}{m_{2}}\right)^{i x }\textrm{d}x, \\
&0<\sigma<1,~~t>0,~~0<\delta_{3}<1,~~0<\delta_{4}<1,~~m_{j}=1,\ldots, [t], \ j=1,2;\notag
\end{align}
with the principal value integral  defined with respect to $x=0$.

It turns out that the asymptotic analysis of the integrals $\tilde{I}_4$ and $\tilde{I}_3$ gives rise respectively to the sum $S_R$ and to the following sum:
\begin{equation}\label{defS3}
S_M(\sigma,t,\delta_2,\delta_3)=\mathop{\sum\sum}_{m_1,m_2\in M(\delta_2,\delta_3,t)} \frac{1}{m_2^{\bar{s}}(m_1+m_2)^s},
\end{equation}
with the set $M$ defined by
\begin{align}\label{1.14}
M(\delta_2,\delta_3,t)=&\left\{ m_j=1,\ldots,[t], \ j=1,2, \ \ \frac{m_2}{m_1}\in \left(\frac{1}{t^{1-\delta_3}-1},t^{1-\delta_2}-1\right) \right\}, \notag \\
&0<\delta_2<1, \qquad 0<\delta_3<1.
\end{align}

\subsection*{An Exact Relation Between $S_R$ and $S_M$}

It turns out that the sums $S_R$ and $S_M$ are related. Indeed, the following exact identities are established in section 3:

\begin{align}\label{1.21}
2\Re&\left\{ 
\sum_{m_{1}=1}^{[t]}\sum_{m_{2}=1}^{[t]}\frac{1}{m_{2}^{\bar{s}} (m_{1}+m_{2})^s} \right\} =S_R(\sigma,t) -\sum_{m=1}^{[t]}\frac{1}{m^{2\sigma }} \notag \\
&\qquad \qquad +2\Re \left\{ \sum_{m=1}^{[t]}\sum_{n=[t]+1}^{[t]+m}\frac{1}{m^{\bar{s} }n^{s}} \right\}, 
\qquad s=\sigma + i t,
\end{align}
and
\begin{equation}\label{1.22}
\sum_{m_{1}=1}^{[t]}\sum_{m_{2}=1}^{[t]}\frac{1}{m_{2}^{\bar{s}} (m_{1}+m_{2})^s} =S_M(\sigma,t,\delta_2,\delta_3) + S_2(\sigma,t,\delta_2) + S_3(\sigma,t,\delta_3),
\end{equation}
 with $S_2$ and $S_3$ defined below;
\begin{align}\label{1.23}
S_2&(\sigma,t,\delta_2)  =
\mathop{\sum\sum}_{(m_{1},m_{2})\in M_2}\frac{1}{m_{2}^{\bar{s}} (m_{1}+m_{2})^s}, \notag \\ &M_2(\delta_2,t)=\Big\{m_{j}=1,\ldots,[t], \ j=1,2, \ \ \frac{m_{2}}{m_{1}} > t^{1-\delta_{2}}\Big\},
\end{align}
and
\begin{align}\label{1.24}
S_3&(\sigma,t,\delta_3)= \mathop{\sum\sum}_{(m_{1},m_{2})\in M_3} \frac{1}{m_{2}^{\bar{s}} (m_{1}+m_{2})^s}\notag \\ 
&M_3(\delta_3,t)=\Big\{m_{j}=1,\ldots,[t], \ j=1,2, \ \ \frac{m_{2}}{m_{1}} < \frac{1}{t^{1-\delta_{3}}-1}\Big\}.
\end{align}

\subsection*{The Rigorous Analysis of $I_1$ and $I_2$}

In section 4 the rigorous analysis of $I_1$ and $I_2$ is presented: it is shown that for $\delta_1$ a sufficiently small positive constant, $I_1$ satisfies the estimate 

\begin{equation} \label{1.8}
I_1(\sigma, t, \delta_1) = \begin{cases}
O \left(  t^{-\sigma+\left(2-\frac{4}{3}\sigma\right)\delta_1} \right), & 0\leq\sigma\leq \frac{1}{2},\\
O \left( t^{-\sigma+\left(\frac{5}{3}-\frac{2}{3}\sigma\right)\delta_1}  \right), & \frac{1}{2} < \sigma< 1,
   \end{cases}, \quad t \to \infty.
\end{equation}

Furthermore, by employing the classical estimates of Atkinson, it is also shown in section 3 that $I_2$ satisfies  the following estimate, for $0<\delta_2<1$:
\begin{equation} \label{1.9}
I_2(\sigma, t,\delta_2) =  \begin{cases}
    O\left( t^{-\sigma + 2 (1 -\sigma )\delta_2} \zeta(2-2\sigma) \right), &  0<\sigma<\frac{1}{2},  \\
     O\left( t^{-\frac{1}{2}+ \delta_2} \ln{t} \right), & \sigma = \frac{1}{2},\\
       O\left( t^{-\sigma + \left(\sigma + \frac{1}{2} \right)\delta_2}\zeta(2\sigma)\right), &  \frac{1}{2}<\sigma< 1,\end{cases} \quad t\to\infty.
\end{equation}

The remaining results presented in this paper are formal.

\subsection*{The Formal Analysis of $\tilde{I}_3$}

The asymptotic analysis as $t\to\infty$ of the integral $J_3$ defined in \eqref{1.11} is discussed in section 5, where it is shown that  $J_3$ involves the following two terms: (i) $J_3^S$, which is due to the contribution of the stationary points, with an error due to the contribution of the lower end point $t^{\delta_2-1}$. \ (ii) $J_3^U$, which is due to the contribution of the upper end point $1-t^{\delta_3-1}$. These terms yield the following representation for $\tilde{I}_3$:

\begin{equation}\label{1.12}
\tilde{I}_3(\sigma,t,\delta_2, \delta_3)=\tilde{I}_3^S(\sigma,t,\delta_2, \delta_3)-\tilde{I}_3^U(\sigma,t,\delta_2, \delta_3),
\end{equation}
where
\begin{align}\label{1.13}
\tilde{I}_3^S(\sigma,t,\delta_2, \delta_3)=2\Re&\left\{ S_M \right\} [1+o(1)] \left[1+O\left(t^{-\delta_{23}}\right)\right], \qquad t \to \infty, \notag\\
& \hspace*{-30mm} s=\sigma+ i  t, \ 0 < \sigma <1, \quad 0 < \delta_2 <1, \  0< \delta_3 <1,  \delta_{23}=\min \{ \delta_2 , \delta_3 \},
\end{align}
and
\begin{align}\label{1.15}
\tilde{I}_3^U&(\sigma,t,\delta_2, \delta_3)=- \sqrt{\frac{2}{\pi}} \Re \Bigg\{ e^{\frac{i\pi}{4}} t^{-\frac{\delta_3}{2}} (1-t^{\delta_3 -1})^{\sigma - \frac{1}{2} +i(t-t^{\delta_3})} t^{i(\delta_3-1)t^{\delta_3}} \notag \\
&\times\mathop{\sum\sum}_{m_1,m_2\in N(\delta_3,t)} \frac{1}{m_1^{s-it^{\delta_{3}}}} \frac{1}{m_2^{\bar{s}+it^{\delta_{3}}}}  \frac{1}{\ln{\left[ \frac{m_2}{m_1} \left( t^{1-\delta_{3}} - 1 \right)  \right]}} [1+o(1)]  
\Bigg\} \notag \\
& \hspace*{10mm} \times \left[  1 + O(t^{-\delta_{23}}) \right], \qquad t \to\infty,\tag{1.24a}
\end{align}
where $N(\delta_3,t)=M_t\cap M_r^c$, with
\begin{equation}
 M_t=\Big\{m_{j}=1,\ldots,[t], \ j=1,2\Big\} \tag{1.24b}
 \end{equation}
and  $M_r^c$ denotes the complement of $M_r$, which is given by
\begin{equation}M_r(\delta_3,t)=\left\{(m_1,m_2), \ \frac{m_{2}}{m_{1}}=t^{\delta_{3}-1}\left(1+ O\left(t^{-\delta_{3}/2}\right)\right)\right\}.\tag{1.24c}
 \end{equation}
\addtocounter{equation}{1}

The set $M_r$ captures the set  of points $(m_1,m_2)$ where the endpoint  of integration $1-t^{\delta_3-1}$ approaches the stationary point. These points  yield 
a contribution which is assumed here to be of lower order (a rigorous analysis is presented in \cite{F}).

\subsection*{The Formal Analysis of $\tilde{I}_4$}

\begin{figure}
\begin{center}
\includegraphics[width=.4\textwidth]{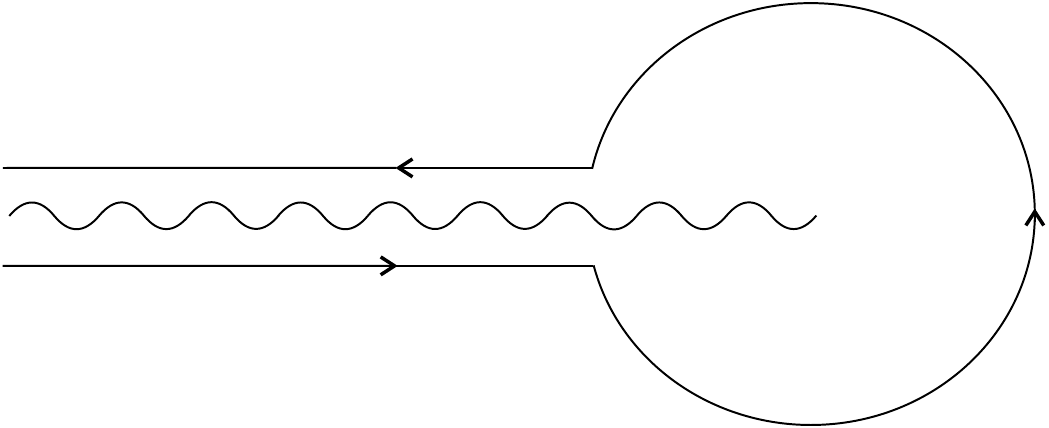}
\end{center}
\caption{The Hankel contour.}
\label{fig-Hankel-1}
\end{figure}

The analysis as $t\to\infty$ of the integral $J_4$ defined in \eqref{1.13a} is discussed in section 6, where using novel analytical and asymptotic techniques implemented in the complex plane, it is shown that $\tilde{I}_4$ also involves two terms: (i) One term can be computed exactly. \ (ii) The second term can be expressed in terms of a sum denoted by $2\Re\left\{S_4^P\right\}$ that can be evaluated exactly in terms of a certain residue computation, and a sum denoted by  $\Re\left\{S_4^{SD}\right\}$ that involves a steepest descent computation. These terms imply the following result for $\tilde{I}_4$:
\begin{align}\label{1.16}
\tilde{I}_{4}&(\sigma, t, \delta_{3}, \delta_{4})=\Bigg[-\sum_{m_1=1}^{[t]}\sum_{m_{2}=1}^{[t]}
\frac{1}{m_{1}^{s}m_{2}^{\bar{s}}} \notag\\
 &+2\Re\left\{S_4^P(\sigma, t, \delta_{3})\right\}+\Re\left\{S_4^{SD}(\sigma, t, \delta_{3})\right\} \Bigg]
\left[1+O(t^{2\delta_{34}-1})\right], \ t\to\infty, \notag\\
& 0<\sigma<1,~0<\delta_{3}<\frac{1}{2},~0<\delta_{4}<\frac{1}{2}, \ \delta_{34}=\min(\delta_3,\delta_4),
\end{align}
where $S_4^P$ is defined by
\begin{align}\label{1.17}
S_4^P(\sigma, t, \delta_{3})=\mathop{\sum\sum}_{m_1,m_2\in M_4(\delta_3,t)}
\frac{1}{m_{1}^{s}m_{2}^{\bar{s}}}e^{-\frac{im_{2}}{m_{1}}t},
\end{align}
with the set $M_4$ defined by
\begin{equation}\label{1.18}
M_4(\delta_3,t)=\Big\{m_{j}=1,\ldots,[t], \ j=1,2, \
~\frac{m_{1}}{m_{2}}\in(t^{1-\delta_{3}},t)\Big\},
\end{equation}
whereas $S_4^{SD}$ is defined by
\begin{align}\label{1.19}
S_4^{SD}(\sigma, t, \delta_{3})=\frac{1}{\pi}\sum_{m_1=1}^{[t]}\sum_{m_{2}=1}^{[t]}\frac{1}{m_{1}^{s}m_{2}^{\bar{s}}}E_{4}^{SD}(t,\delta_{3},M),
\end{align}
with $E_{4}^{SD}(t,\delta_{3},M)$ defined by
\begin{equation}\label{1.20}
E_{4}^{SD}=\int_{H_{1}}\frac{e^{t^{\delta_{3}}[\omega-\frac{\pi}{2}+i\ln(Mt^{\delta_{3}}\omega)] }      }{\omega[\frac{\pi}{2}-i\ln(Mt^{\delta_{3}}\omega)]}d\omega , \quad M=\frac{m_1}{m_2 t},
\end{equation}
and $H_1$ denoting the Hankel contour, with a branch cut along the negative real axis, see figure \ref{fig-Hankel-1},
\begin{align*}
H_{1}=\left\{r e^{- i \pi}|1<r<\infty\right\}\cup\left\{e^{ i \theta}|-\pi<\theta<\pi\right\}\cup\left\{r e^{ i \pi}|1<r<\infty\right\}.
\end{align*}

\subsection*{Formal Asymptotics for a certain Riemann-type Exponential Sum}

It turns out that the sums $S_3$ and $S_4^P$ defined by \eqref{1.24} and \eqref{1.17} are related: 
the following relation and estimate are derived in \cite{F} and \cite{KF}, respectively:
\begin{equation}\label{1.25}
S_3(\sigma,t,\delta_3)= S_4^P(\sigma, t, \delta_{3}) \left[1+ O\left(t^{2\delta_3-1} \right)\right],\qquad t\to \infty,
\end{equation}
and
\begin{align} \label{1.26}
 \sum_{m=1}^{[t]}\sum_{n=[t]+1}^{[t]+m}\frac{1}{m^{\bar{s}}n^{s} } 
= \begin{cases}
O \left(  t^{\frac{1}{2}-\frac{5}{3}\sigma} \ln{t} \right), & 0\leq\sigma\leq \frac{1}{2},\\
O \left(  t^{\frac{1}{3}-\frac{4}{3}\sigma} \ln{t}  \right), & \frac{1}{2} < \sigma< 1,
\end{cases} \qquad t\to \infty.
\end{align}

We restrict our consideration to the case $\sigma=\frac{1}{2}$, with $\delta_1, \delta_3,  \delta_4$ arbitrarily small positive constants and $0<\delta_2<1$. 

Adding  the expressions for $\tilde{I}_3$ and $\tilde{I}_4$ given by equations \eqref{1.12} and \eqref{1.16}, and then using equations \eqref{1.21}, \eqref{1.22} and \eqref{1.25} to simplify the resulting expression, we find 
\begin{align}\label{1.27}
&\tilde{I}_3(\sigma,t,\delta_2, \delta_3) + \tilde{I}_4(\sigma,t,\delta_3, \delta_4) \sim \notag \\
& \Re\left\{S_4^{SD}(\sigma, t, \delta_{3})\right\} -\tilde{I}_3^U(\sigma,t, \delta_3) - 
2\Re \left\{S_2(\sigma,t,\delta_2) \right\} -\sum_{m=1}^{[t]}\frac{1}{m^{2\sigma }}, \quad t\to\infty, \notag\\
&\sigma=\frac{1}{2}, \ 0<\delta_2<1, \ \delta_3, \ \delta_4 \text{ small positive constants} .
\end{align}

Comparing this equation with  \eqref{1.5}, namely
\begin{equation}\label{1.28}
\tilde{I}_3+\tilde{I}_4 \sim - I_2 - \mathcal{G}(\sigma,t), \qquad t\to\infty,
\end{equation}
evaluated at $\sigma=\frac{1}{2}$, where the asymptotic behaviour of $\mathcal{G}(\sigma,t)$ as $t\to \infty$  is given by \eqref{1.4b}, we find the following: first, the leading asymptotics of $\Re\left\{S_4^{SD}\right\}$ is given by the leading asymptotics of $\tilde{I}_3^U$ (as noted in remark 6.1, this result can be verified explicitly using the steepest descent computation). Second,
\begin{equation}\label{1.29}
\Re \left\{S_2\left(\frac{1}{2},t,\delta_2\right) \right\} = I_2\left(\frac{1}{2},t, \delta_2\right) +O(1), \qquad t \to \infty.
\end{equation}
Making the change of variables $m_1+m_2=n$ and $m_2=m$ in the definition of $S_2$ given in \eqref{1.23}, and  using the estimate for $I_2$ given in \eqref{1.9}  evaluated at $\sigma=\frac{1}{2}$, equation \eqref{1.29} yields
\begin{equation}\label{1.30}
\sum_{m=\left[t^{1-\delta _{2}}\right]}^{[t]}\sum_{n=m+1}^{\left[m\left(1+t^{\delta _{2}-1}\right)\right]}\frac{1}{n^{\frac{1}{2}+i t} m^{\frac{1}{2}-it}} =  O\left( t^{-\frac{1}{2}+ \delta_2} \ln{t} \right) +O(1), \quad t\to\infty.
\end{equation}
For the particular case of $\delta_2=1/2+\epsilon$,  the above equation becomes 
\begin{equation}\label{R1}
\sum_{m=\left[t^{\frac{1}{2}-\epsilon}\right]}^{[t]}\sum_{n=m+1}^{\left[m\left(1+t^{\epsilon-\frac{1}{2}}\right)\right]}\frac{1}{n^{\frac{1}{2}+i t} m^{\frac{1}{2}-it}}  = O\left( t^{\epsilon} \ln{t} \right), \quad t \to \infty, \ \ \epsilon>0.
\end{equation}

This equation is the analogue of the Lindel\"of hypothesis for the Riemann-type double sum defined by the lhs of the above equation.

\subsection*{Further Developments}

The combination of the Atkinson asymptotic result and the asymptotic analysis of various integrals employed above suggests further development. Some of these developments are briefly discussed in section 7: the analysis of carefully chosen integrals that involve $|\zeta(1/2+i\tau t)|^2$ yields asymptotic estimates for a variety of double exponential sums. In this analysis, on the one hand, $|\zeta(1/2+i\tau t)|^2$ is replaced by its leading asymptotic double sum and an evaluation of the resulting integral is performed;  and on the other hand, the original integral is estimated using Atkinson's asymptotic result employed for the estimation of $I_2$. In section 7 only a couple of such integrals are formally analysed; a rigorous thorough investigation will be presented in \cite{F}.

The above carefully chosen integrals also provide a promising avenue of an  investigation first suggested in an earlier unpublished preprint of Kalimeris, Lenells and the author.  Unfortunately, the integral analysed in this preprint was incorrect due to the fact that the earlier asymptotic analysis of $I_4$ missed the important term $E_4$. Hopefully, by a judicious choice of an integral of the type analysed in section 7, this approach can be implemented: consider an integral involving $|\zeta(1/2+iy)|^2$ and a combination of gamma functions that has two properties; first, the oscillations of the above combination `captures' the oscillations of $|\zeta(1/2+iy)|^2$, and second the modulus of this combination is `sufficiently small'. This yields an equation similar to \eqref{R1};  the first property yields a double exponential sum on the lhs of this equation that is a slight variant of $|\zeta(1/2+iy)|^2$ and the second property yields an estimate on the rhs of this equation that improves the best current estimate for the asymptotics of the Riemann zeta function. 

\section{A Singular Integral Equation for $\left|\zeta(s)\right|^2$ and its asymptotic form} \label{sec2}
In this section we derive equations \eqref{1.3}  and \eqref{1.5}.

\begin{proposition} Let $u=\sigma+it$, $\sigma>1$, $\omega > -1$. The following identity is valid:
\begin{equation}\label{2.1}
\Gamma(u)(1+\omega)^{-u}=\frac{1}{2i\pi}\int_{c-i\infty}^{c+i\infty}
\Gamma(u+z)\Gamma(-z)\omega^{z}dz,~~-\sigma<c<0.
\end{equation}
\end{proposition}

\textbf{Proof}
We will prove (\ref{2.1}) for $\omega\in (-1, 1)$; analytic continuation then implies that it is valid for all $\omega > -1$.
Consider a clockwise contour $C$ enclosing $0,1,2,\ldots$,
but not $-\sigma,-\sigma - 1,\ldots$. Then $\Gamma(u+z)$ does not have
any poles in the domain enclosed by $C$, whereas $\Gamma(-z)$ has poles
at $z=n$ for $n=0,1,2,\ldots$, with residues given by
\begin{equation*}
\mbox{Res}_{n} \Gamma(-z)=\frac{(-1)^{n+1}}{n!},~~n=0,1,\ldots.
\end{equation*}
Hence
$$\frac{1}{2i\pi}\int_{C}\Gamma(u+z)\Gamma(-z)\omega^{z}dz=
\sum_{n=0}^{\infty}\frac{(-1)^{n}}{n!}\Gamma(u+n)\omega^{n}=
(1+\omega)^{-u}\Gamma(u),
$$
where in the last step we have used the binomial theorem and the
fact that $|\omega|<1$.

By Stirling's approximation, we have the formula
$$\Gamma(z)=e^{-z}z^{z}\left(\frac{2\pi}{z}\right)^{\frac{1}{2}}(1+O(z^{-1})),~~
z\rightarrow\infty,~~|\mbox{arg}z|\leq \pi-\delta.
$$
Thus,
$$\frac{\Gamma(z+a)}{\Gamma(z+b)}\sim e^{b-a} z^{a-b},~~z\rightarrow\infty.
$$
This estimate, together with the identity
$$\Gamma(-z)=-\frac{\pi}{\sin(\pi z)\Gamma(1+z)},
$$
yield
$$ \Gamma(-z)\Gamma(z+u)\sim z^{u-1}\frac{\pi}{2i(e^{-i\pi z}-e^{i\pi z})}.
$$
Thus, the product of the above gamma functions decays, and since $\omega^{z}=e^{z\ln \omega}$ also decays for $|\omega|<1$, it follows
that we can replace the contour $C$ by the straight line from $c-i\infty$
to $c+i\infty$,  and then replace the condition $|\omega|<1$ by the condition $\omega>-1$.   \textbf{QED}

%\begin{corollary}
%Multiplying \eqref{2.1} by $\alpha^{-u}/\Gamma(u)$, $\alpha>0$, $u=\sigma+it$,
%$\sigma>1$, $t>0$, letting $\omega=m/\alpha$ and $-\sigma<c<-1$, and
%summing from $m=1$ to $m=\infty$, we find
%\begin{align}\label{2.2}
%&\zeta_{1}(u,\alpha)=\frac{1}{2i\pi}\int_{c-i\infty}^{c+i\infty}
%\frac{\Gamma(u+z)\Gamma(-z)}{\Gamma(u)}\zeta(-z)\alpha^{-u-z}dz,~~
%-\sigma<c<-1,\notag\\
%&\sigma>1,~~t>0,
%\end{align}
%where $\zeta_{1}(u,\alpha)$, $u\in\mathbb{C}$, $\alpha\geq 0$,
%denotes 

%Deforming the above contour past the pole at $z=-1$, equation \eqref{2.2}
%becomes
%\begin{align}\label{2.4}
%&\zeta_{1}(u,\alpha)=\frac{1}{2i\pi}\int_{c-i\infty}^{c+i\infty}
%\frac{\Gamma(u+z)\Gamma(-z)}{\Gamma(u)}\zeta(-z)\alpha^{-u-z}dz
%+\frac{\alpha^{1-u}}{u-1},~~
%-1<c<0,\notag\\
%&\sigma>1, ~~ t>0.
%\end{align}
%\end{corollary}

\begin{proposition}
Let $s=\sigma+it, \ \sigma>1, \ t\geq 0$. Define the modified Hurwitz function by the expression
\begin{equation}\label{2.2}
\zeta_{1}(s,\alpha)=\sum_{n=1}^{\infty}\frac{1}{(n+\alpha)^s},~~ \Re s>1,
~~\alpha\geq 0.
\end{equation}
This function which is related to the usual Hurwitz function $\zeta(s,\alpha)$ via the relation $$\zeta(s,\alpha)=\frac{1}{\alpha}+\zeta_{1}(s,\alpha),$$ can be defined
  via analytic continuation for all $s\in \mathbb{C}$. The following identity is valid:
\begin{align}\label{2.5}
&\left| \zeta_{1}(s,\alpha)\right|^2
=\zeta_{1}(2\sigma,\alpha)+\frac{1}{2i\pi}\int_{c-i\infty}^{c+i\infty}
\left(\frac{\Gamma(s+z)}{\Gamma(s)}+\frac{\Gamma(\bar{s}+z)}{\Gamma(\bar{s})}\right)
\notag\\
&\times \Gamma(-z)\zeta(-z)\zeta_{1}(2\sigma+z,\alpha)dz,
\end{align}
where $\bar{s}=\sigma-it$ and 
$$-\sigma <c<-1,~~\sigma>1,~~t>0.$$

\end{proposition}

\textbf{Proof}  
Define the function $f(u,v,\alpha)$ by 
\begin{equation}\label{2.4n}
f(u,v,\alpha)=\sum_{n=1}^{\infty} \sum_{m=1}^\infty (n+\alpha)^{-u} (n+m+\alpha)^{-v},  \qquad \alpha\geq 0, \ \Re u>1, \ \Re v>1.
\end{equation}

Letting in \eqref{2.1}  $\omega=\frac{m}{n+\alpha}$, we find 
\begin{align*}
\Gamma(u)(n+\alpha)^u (m+ n +\alpha)^{-u}=\frac{1}{2i\pi}\int_{c-i\infty}^{c+i\infty}
\Gamma(u+z)\Gamma(-z) m^z (n+\alpha)^{-z} dz,\\
-\sigma<c<0.
\end{align*}
Multiplying this equation by $\Gamma(u)^{-1} (n+\alpha)^{-u} (n+\alpha)^{-v}$, we find
\begin{equation}\label{2.5n}
(n+\alpha)^{-v} (m+ n +\alpha)^{-u}=\frac{1}{2i\pi}\int_{c-i\infty}^{c+i\infty}
\frac{\Gamma(u+z)}{\Gamma(u)}\Gamma(-z) m^z (n+\alpha)^{-(z+u+v)} dz.
\end{equation}
Summing over $m, \ n$ and using the definition of $f$ we obtain 
\begin{equation}\label{2.6n}
f(v,u,\alpha)=\frac{1}{2i\pi}\int_{c-i\infty}^{c+i\infty}
\frac{\Gamma(u+z)}{\Gamma(u)}\Gamma(-z) \zeta(-z) \zeta_1(u+v+z,\alpha) dz.
\end{equation}
Using straightforward calculations, the definitions of $\zeta_1(u,\alpha)$ and $f(u,v,\alpha)$ yield the identity
\begin{equation}\label{2.7n}
\zeta_1(u,\alpha)\zeta_1(v,\alpha)=\zeta_1(u+v,\alpha)+f(u,v,\alpha)+f(v,u,\alpha).
\end{equation}
 Indeed, for $\Re u>1$ and $\Re v>1$, we have
\begin{align*}
\zeta_1(u,\alpha) \zeta_1(v, \alpha) 
= &\; \sum_{n = 1}^\infty \sum_{m = 1}^\infty (n + \alpha)^{-u}  (m + \alpha)^{-v} 
=  \sum_{n = 1}^\infty \sum_{m = 1}^{n-1} (n + \alpha)^{-u}  (m + \alpha)^{-v} 
	\\
& + \sum_{n = 1}^\infty (n + \alpha)^{-u-v}
+ \sum_{n = 1}^\infty \sum_{m = n+1}^{\infty} (n + \alpha)^{-u}  (m + \alpha)^{-v} .
\end{align*}
The identity (\ref{2.7n}) follows because of the following relations:
\begin{align*}
\sum_{n = 1}^\infty \sum_{m = 1}^{n-1} (n + \alpha)^{-u}  (m + \alpha)^{-v} 
& = \sum_{m = 1}^\infty \sum_{n = m+1}^\infty (n + \alpha)^{-u}  (m + \alpha)^{-v} 
	\\
& = \sum_{m = 1}^\infty \sum_{n = 1}^\infty (m + \alpha)^{-v}  (n + m + \alpha)^{-u}  
= f(v,u,\alpha).
	\\
 \sum_{n = 1}^\infty \sum_{m = n+1}^{\infty} (n + \alpha)^{-u}  (m + \alpha)^{-v}
&= \sum_{n = 1}^\infty \sum_{m = 1}^\infty (n + \alpha)^{-u}  (n + m + \alpha)^{-v}
= f(u,v,\alpha),
\end{align*}
and
$$\sum_{n = 1}^\infty (n + \alpha)^{-u-v} = \zeta_1(u+v, \alpha).$$

Replacing in \eqref{2.7n} the functions $f$ via the representations \eqref{2.6n}, and then replacing in the resulting equation $u$ and $v$ by $s$ and $\bar{s}$ we obtain equation \eqref{2.5}.
 \textbf{QED}

\begin{lemma} \label{l2.1}
The Riemann zeta function $\zeta(s)$ satisfies the identity
\begin{align}
\left|\zeta(s)\right|^2 &= \mathcal{G}(\sigma,t) \nonumber\\
&+ \frac{1}{2i\pi} \int_{c-i\infty}^{c+i\infty} \left( \frac{\Gamma(s+z)}{\Gamma(s)} + \frac{\Gamma(\bar{s}+z)}{\Gamma(\bar{s})} \right) \Gamma(-z) \zeta(-z) \zeta(2\sigma+z)\textrm{d}z, \nonumber
\end{align}
\begin{equation} \label{2.3}
s = \sigma +it, \quad \sigma>0, \quad t>0, \quad \max{(-1,-\sigma)} < c < \min{(0,1-2\sigma)},
%\end{gathered}
\end{equation}
where the function $\mathcal{G}(\sigma,t)$ is defined by the following expressions:
\begin{align} \label{2.4}
\mathcal{G}(\sigma,t) = \zeta(2\sigma) &+ \left( \frac{\Gamma(1-\bar{s})}{\Gamma(s)} + \frac{\Gamma(1-s)}{\Gamma(\bar{s})} \right) \Gamma(2\sigma-1) \zeta(2\sigma-1) \nonumber  \\
&+ \frac{2(\sigma -1)\zeta(2\sigma-1)}{(\sigma - 1)^2 + t^2}, \qquad \sigma \ne \frac{1}{2},
\end{align}

\begin{equation} \label{2.5b}
\mathcal{G}\left(\frac{1}{2},t\right) = \Re\left(\Psi\left(\frac{1}{2} + it\right)\right) + 2\gamma - \ln{2\pi} + \frac{2}{1+4t^2},
\end{equation}
with $\Psi(z)$ denoting the digamma function, i.e.
\begin{equation}  \label{2.6}
\Psi(z) = \frac{\frac{\mathrm{d}}{\mathrm{dz}}\Gamma(z)}{\Gamma(z)}, \quad z \in \mathbb{C},
\end{equation}
and $\gamma$ denoting the Euler constant.
\end{lemma}
\textbf{Proof}
Let us first recall that the function $\mu(\sigma)$ defined as the lower bound of numbers $\xi$ such that $\zeta(\sigma + it) = O(|t|^\xi)$ as $t�\to \pm \infty$ satisfies (see Chapter V of Titchmarsh)
$$\mu(\sigma) = \begin{cases} 0, & \sigma > 1, \\
\frac{1}{2} - \sigma, & \sigma < 0,
\end{cases} $$
and $\mu(\sigma) \leq \frac{1-\sigma}{2}$ for $\sigma \in (0,1)$. This provides basic bounds for $\zeta(s)$ which can be used to justify the contour deformations below.

Letting $\alpha=0$ in equation \eqref{2.5} we find
\begin{equation} \label{2.7}
%\begin{gathered}
\left|\zeta(s)\right|^2 = \zeta(2\sigma)   \nonumber
\end{equation}
\begin{equation}
+ \frac{1}{2i\pi} \int_{c-i\infty}^{c+i\infty} \left( \frac{\Gamma(s+z)}{\Gamma(s)} + \frac{\Gamma(\bar{s}+z)}{\Gamma(\bar{s})} \right) \Gamma(-z) \zeta(-z) \zeta(2\sigma+z)\textrm{d}z, \nonumber
\end{equation}
\begin{equation}
s = \sigma +it, \quad \sigma>1, \quad t>0, \quad -\sigma < c < -1.
%\end{gathered}
\end{equation}
The integrand in (\ref{2.7}) has simple poles at $z=1-2\sigma$ and $z = -1$ (because $\zeta(s)$ has a simple pole at $s = 1$) and at $z = -\sigma \pm it$ (because 
$\Gamma(s)$ has a simple pole at $s = 0$). Let us first consider (\ref{2.7}) with $\sigma$ slightly larger than $1$. Then $1-2\sigma < -\sigma$, so the three poles $z=1-2\sigma$, $z = -s$, and $z = -\bar{s}$ lie to the left of the integration contour, whereas the pole $z = -1$ lies to the right of the contour. 
We begin by deforming the contour to the right by increasing $c$ to a value slightly larger than $-1$. During this deformation, a residue contribution (called $R_2$ below) is picked up from the pole at $z = -1$. We next decrease the value of $\sigma$. For $\sigma = 1$, the three poles $1-2\sigma$ and $-\sigma \pm it$ all lie on the vertical line with real part $-1$. As we keep decreasing $\sigma$, the pole at $z =1-2\sigma$ lies to the right of the poles $-\sigma \pm it$ and will eventually (for $\sigma = (1-c)/2$) pass through the contour, which generates another residue contribution (called $R_1$ below). In this final configuration, the poles $-1$ and $-\sigma \pm it$ lie to the left of the contour, while $1-2\sigma$ lies to the right of the contour.    It should be noted that because of the reality property of the zeta-function, i.e. $\bar{\zeta}(s)=\zeta(\bar{s})$, we have $$|\zeta(s)|^2=\zeta(\sigma+it)\zeta(\sigma-it).$$ 
Thus, both sides of \eqref{2.3} and \eqref{2.7} admit a meromorphic continuation to the complex $\sigma$-plane and hence we indeed can move $-\sigma$ to the right after we have deformed the contour in \eqref{2.7} to the contour with $c<-1$. Analytic continuation leads to the following formula, valid for $0<\sigma<1,\  \sigma<1/2 \  -\sigma<c<1-2\sigma$:
\begin{equation} \label{2.8}
%\begin{gathered}
\left|\zeta(s)\right|^2 = \zeta(2\sigma) + R_1 +R_2   \nonumber
\end{equation}
\begin{equation}
+ \frac{1}{2i\pi} \int_{c-i\infty}^{c+i\infty} \left( \frac{\Gamma(s+z)}{\Gamma(s)} + \frac{\Gamma(\bar{s}+z)}{\Gamma(\bar{s})} \right) \Gamma(-z) \zeta(-z) \zeta(2\sigma+z)\textrm{d}z, \nonumber
\end{equation}
\begin{equation}
s = \sigma +it, \ \ \sigma \in (0,1), \ \ \sigma \neq \frac{1}{2}, \ \ t>0, \ \  -\sigma < c < 1-2\sigma ,
%\end{gathered}
\end{equation}
where $R_1$ and $R_2$ are defined by the following expressions:
\begin{equation*}
R_1 = \underset{z=1-2\sigma}{\operatorname{Res}} {\left( \frac{\Gamma(s+z)}{\Gamma(s)} + \frac{\Gamma(\bar{s}+z)}{\Gamma(\bar{s})} \right) \Gamma(-z) \zeta(-z) \zeta(2\sigma+z)},
\end{equation*}
\begin{equation*}
R_2 = - \underset{z=-1}{\operatorname{Res}} {\left( \frac{\Gamma(s+z)}{\Gamma(s)} + \frac{\Gamma(\bar{s}+z)}{\Gamma(\bar{s})} \right) \Gamma(-z) \zeta(-z) \zeta(2\sigma+z)}.
\end{equation*}
The above residues can be computed as follows:
\begin{equation}\label{2.9}
R_1 =  \left( \frac{\Gamma(1-\bar{s})}{\Gamma(s)} + \frac{\Gamma(1-s)}{\Gamma(\bar{s})} \right) \Gamma(2\sigma - 1) \zeta(2\sigma - 1) \underset{z=1-2\sigma}{\operatorname{Res}} \zeta(2\sigma+z),
\end{equation}
\begin{equation}\label{2.10}
R_2 = -  \left( \frac{\Gamma(s-1)}{\Gamma(s)} + \frac{\Gamma(\bar{s}-1)}{\Gamma(\bar{s})} \right) \Gamma(1) \zeta(2\sigma - 1) \underset{z=-1}{\operatorname{Res}} \zeta(-z).
\end{equation}
Using the identities
\begin{equation} \label{2.11}
\Gamma(u) = (u-1)\Gamma(u-1),~u\in\mathbb{C}; \quad \Gamma(1)=1,
\end{equation}
and
\begin{equation} \label{2.12}
\zeta(1+\varepsilon) = \frac{1}{\varepsilon} + \gamma + O(\varepsilon), \quad \varepsilon\to 0,
\end{equation}
equations \eqref{2.9} and \eqref{2.10} become
\begin{equation}\label{2.13}
R_1 =  \left( \frac{\Gamma(1-\bar{s})}{\Gamma(s)} + \frac{\Gamma(1-s)}{\Gamma(\bar{s})} \right) \Gamma(2\sigma - 1) \zeta(2\sigma - 1),
\end{equation}
\begin{equation}\label{2.14}
R_2 = \frac{2(\sigma-1)}{(\sigma-1)^2+t^2} ~ \zeta(2\sigma - 1).
\end{equation}
Substituting in \eqref{2.8} the expressions for $R_1$ and $R_2$ given by equations \eqref{2.13} and \eqref{2.14} we find equation \eqref{2.4}.

We next compute the limit of the rhs of equation \eqref{2.4} as $\sigma\to\frac{1}{2}$. For this purpose, in addition to equation \eqref{2.12} we will also use the following identities:
\addtocounter{equation}{1}
\begin{equation}\label{2.15}
\Gamma(\varepsilon) = \frac{1}{\varepsilon} - \gamma + O(\varepsilon), \quad \varepsilon\to 0,\tag{\theequation a}
\end{equation}
and
\begin{equation}\label{2.16}
\zeta(\varepsilon) = -\frac{1}{2} \left(1+\varepsilon\ln{2\pi} + O(\varepsilon^2)\right), \quad \varepsilon\to 0. \tag{\theequation b}
\end{equation}
Letting
$
\sigma = \frac{1}{2} + \frac{\varepsilon}{2},
$
we find
\begin{equation*}
 \frac{\Gamma(1-\bar{s})}{\Gamma(s)} =  \frac{\Gamma(1-\sigma+it)}{\Gamma(\sigma+it)} = \frac{\Gamma\left(\frac{1}{2}+it -\frac{\varepsilon}{2}\right)}{\Gamma\left(\frac{1}{2}+it + \frac{\varepsilon}{2}\right)}   \nonumber
\end{equation*}
\begin{equation*}
= \frac{1 - \frac{\varepsilon}{2}\frac{\Gamma'}{\Gamma}\left(\frac{1}{2}+it\right) + O(\varepsilon^2)}{1
+ \frac{\varepsilon}{2}\frac{\Gamma'}{\Gamma}\left(\frac{1}{2}+it\right) + O(\varepsilon^2)}.
\end{equation*}
Thus,
\begin{equation}  \label{2.17}
 \frac{\Gamma(1-\bar{s})}{\Gamma(s)} = 1 - \varepsilon \Psi\left(\frac{1}{2}+it\right) + O(\varepsilon^2).
\end{equation}
Equations (2.20) imply
\begin{equation}  \label{2.18}
 \Gamma(2\sigma - 1)\zeta(2\sigma - 1) = \Gamma(\varepsilon)\zeta(\varepsilon) = - \frac{1}{2} \left(\frac{1}{\varepsilon} + \ln{2\pi} - \gamma + O(\varepsilon) \right),  \nonumber
\end{equation}
\begin{equation}
\quad \varepsilon \to 0.
\end{equation}
Using in equation \eqref{2.4}, equations \eqref{2.12} and (2.20) we find
\begin{equation*}
\mathcal{G}\left(\frac{1}{2}+\frac{\varepsilon}{2},t\right) = \frac{1}{\varepsilon} + \gamma + \frac{2}{1+4t^2}+ O(\varepsilon)  \nonumber
\end{equation*}
\begin{equation*}
 - \left[ 1 - \varepsilon \Re{\left(\Psi\left(\frac{1}{2}+it\right)\right)} + O(\varepsilon^2) \right] \left[ \frac{1}{\varepsilon} + \ln{2\pi} -\gamma + O(\varepsilon) \right],~~\varepsilon \to 0.
\end{equation*}
Hence,
\begin{equation*}
\mathcal{G}\left(\frac{1}{2}+\frac{\varepsilon}{2},t\right) = \Re\left(\Psi\left(\frac{1}{2} + it\right)\right) + 2\gamma - \ln{2\pi} + \frac{2}{1+4t^2}+ O(\varepsilon), ~~\varepsilon \to 0,
\end{equation*}
and equation \eqref{2.5b} follows.
\textbf{QED}

\vphantom{a}

%For the derivation of the large $t$ asymptotics of the function $\mathcal{G}$ defined in \eqref{2.4} and \eqref{2.5} the following remark is useful.

%Using equation (2.8) and employing the Plemelj formulae, the following basic equation can be derived.

Equation \eqref{1.3} follows from (2.8) with the aid of the so-called Plemelj formulas: consider a closed curve $L$ dividing the complex $z$-plane in the domain $D^+$ inside $L$ and the domain $D^-$ outside $L$; let $f(\tau), \ \tau\in L$ be a H\"older function. Then,
$$\lim_{z\to\xi\in L} \frac{1}{2i\pi}\int_L \frac{f(\tau)}{\tau-z}d\tau = \pm \frac{f(\xi)}{2} + \frac{1}{2i\pi} \oint_L \frac{f(\tau)}{\tau-\xi}d\tau , \qquad z\in D^\pm,$$
where the principal value integral is with respect to $\tau=\xi$.
\begin{theorem} \label{t2.1}
The Riemann zeta function $\zeta(s)$ satisfies the integral equation
\begin{equation} \label{2.23}
\frac{t}{\pi} \oint_{-\infty}^{\infty} \Re{\left\{ \frac{\Gamma(it-i\tau t)}{\Gamma(\sigma + it)}\Gamma(\sigma + i\tau t) \right\}} \left| \zeta(\sigma+i\tau t)\right|^2 \textrm{d}\tau + \mathcal{G}(\sigma,t) = 0,   \nonumber
\end{equation}
\begin{equation}
\quad 0<\sigma<1, \quad t>0,
\end{equation}
where $\mathcal{G}(\sigma,t)$ is defined by equations \eqref{2.4} and \eqref{2.5b}, and the principal value integral in equation \eqref{2.23} is defined with respect to $\tau=1$. %i.e.,
%\begin{equation*}
%\oint_{-\infty}^{\infty}\textrm{d}\tau = \lim_{\varepsilon\to 0} \left( \int_{-\infty}^{1-\varepsilon}\textrm{d}\tau + \int_{1+\varepsilon}^{\infty}\textrm{d}\tau \right).
%\end{equation*}
\end{theorem}

\textbf{Proof}
We take the limit of equation (2.8) as $c \downarrow -\sigma$. In this limit the poles $z=-s$ and $z=-\bar{s}$ approach the contour of integration from the left. Using  Plemelj's formulae we find the following equation:
\begin{equation} \label{2.24}
\left|\zeta(s)\right|^2 = \mathcal{G}(\sigma,t) + P_1+P_2  \nonumber
\end{equation}
\begin{equation}
+ \frac{1}{2i\pi} \oint_{-\sigma-i\infty}^{-\sigma+i\infty} \left( \frac{\Gamma(s+z)}{\Gamma(s)} + \frac{\Gamma(\bar{s}+z)}{\Gamma(\bar{s})} \right) \Gamma(-z) \zeta(-z) \zeta(2\sigma+z)\textrm{d}z,
\end{equation}
where
\begin{equation}  \label{2.25}
P_1 = \frac{i\pi}{2i\pi} \underset{z=-s}{\operatorname{Res}} {\left( \frac{\Gamma(s+z)}{\Gamma(s)} + \frac{\Gamma(\bar{s}+z)}{\Gamma(\bar{s})} \right) \Gamma(-z) \zeta(-z) \zeta(2\sigma+z)},
\end{equation}
\begin{equation}  \label{2.26}
P_2 = \frac{i\pi}{2i\pi} \underset{z=-\bar{s}}{\operatorname{Res}} {\left( \frac{\Gamma(s+z)}{\Gamma(s)} + \frac{\Gamma(\bar{s}+z)}{\Gamma(\bar{s})} \right) \Gamma(-z) \zeta(-z) \zeta(2\sigma+z)}.
\end{equation}
Employing equation \eqref{2.15} it follows that
\begin{equation}
P_1=P_2=\frac{1}{2}|\zeta(s)|^2.
\end{equation}
Thus, equation \eqref{2.24} simplifies to the equation
\begin{align} \label{2.28}
&\frac{1}{2i\pi}  \oint_{-\sigma-i\infty}^{-\sigma+i\infty} \left( \frac{\Gamma(s+z)}{\Gamma(s)} + \frac{\Gamma(\bar{s}+z)}{\Gamma(\bar{s})} \right) \Gamma(-z) \zeta(-z) \zeta(2\sigma+z)\textrm{d}z+ \mathcal{G}(\sigma,t) =0,\notag\\
&s=\sigma+it,~~t>0,~~0<\sigma<1.
\end{align}
Letting
$
z = - \sigma + i\tau t,
$
equation \eqref{2.28} becomes
\begin{equation} \label{2.29}
\frac{t}{2\pi} \oint_{-\infty}^{\infty}\left( \Gamma(it+i\tau t) \frac{\Gamma(\sigma-i\tau t)}{\Gamma(\sigma+i t)} + \Gamma(-it+i\tau t) \frac{\Gamma(\sigma-i\tau t)}{\Gamma(\sigma-i t)}\right) \left| \zeta(\sigma+i\tau t)\right|^2 \textrm{d}\tau  \nonumber
\end{equation}
\begin{equation}
+ \mathcal{G}(\sigma,t) = 0,
\end{equation}
where now the principal value integral is defined with respect to $\tau=1$ and $\tau=-1$.

Letting $\tau \to -\tau$ in the first part of the integral in the lhs of equation \eqref{2.29} we find
\begin{equation*}
\frac{t}{2\pi} \oint_{-\infty}^{\infty} \Gamma (it+i\tau t) \frac{\Gamma (\sigma-i\tau t)}{\Gamma (\sigma+i t)} \left| \zeta(\sigma+i\tau t)\right|^2 \textrm{d}\tau   \nonumber
\end{equation*}
\begin{equation*}
= \frac{t}{2\pi} \oint_{-\infty}^{\infty} K(\sigma, t, \tau ) \left| \zeta(\sigma+i\tau t)\right|^2 \textrm{d}\tau,
\end{equation*}
where the kernel $K$ is defined by
\begin{equation} \label{2.30}
K(\sigma, t, \tau) =  \Gamma(it-i\tau t) \frac{\Gamma(\sigma+i\tau t)}{\Gamma(\sigma+i t)}, \quad 0<\sigma< 1, \quad t>0,~~ \tau \in (-\infty,\infty).
\end{equation}
Hence, equation \eqref{2.29} can be rewritten in the form of equation \eqref{2.23}.

\textbf{QED}

\vphantom{a}
  \begin{remark}
  Due to  the fact that the integrand of the integral in the rhs of \eqref{2.3} is analytic, it is possible to derive equation \eqref{2.23} without the use of the Plemelj formulae: replace the vertical line through $z=c$ with the contour that starts at $-\sigma-i\infty$, goes over a small right semicircle of radius $\epsilon$ and centred at the pole $\zeta=-\sigma-it$, passes through $-\sigma$, goes over a small right semicircle of radius $\epsilon$ and centred at the pole $\zeta=-\sigma+it$, and ends up at $-\sigma+i\infty$. Then, taking the limit $\epsilon\to 0$ we find \eqref{2.24}.
\end{remark}

\begin{remark}
Equation \eqref{2.15} implies the following estimate for the singularity of $\Gamma(it-i\tau t)$ at $\tau = 1$:
\begin{equation} \label{2.31}
K(\sigma, t, \tau) =   -\frac{1}{i t (\tau - 1)} + O(1), \quad \tau\to 1.
\end{equation}
The integral equation \eqref{2.23} involves the {\it real} part of $K$ which according to equation \eqref{2.31} is non-singular at $\tau=1$, hence {\it no} principal value integral is needed in \eqref{2.23}. However, in the analysis that follows, instead of the real part of $K$ we will first compute $K$, and for this reason it is useful to retain the principal value.
\end{remark}

%\section{A Simplification of the Integral Equation for $\left|\zeta(s)\right|^2$} \label{sec3}

In order to analyse the large $t$ behaviour of equation \eqref{2.23} we first
use the fact that for $-\infty<\tau<+\infty$ the gamma functions occurring
in the lhs of (2.23) decay exponentially, unless $-t^{\delta_{1}-1}<\tau<1+t^{\delta_{4}-1}$, where $\delta_1$ and $\delta_4$ are positive constants.

\begin{lemma} \label{l3.1}
The Riemann zeta function $\zeta(s)$ satisfies the integral equation
\begin{align} \label{3.1}
&\frac{t}{\pi} \oint_{-t^{\delta_1-1}}^{1+t^{\delta_4-1}} \Re {\left\{ \frac{\Gamma(it-i\tau t)}{\Gamma(\sigma + it)}\Gamma(\sigma + i\tau t) \right\}} \left| \zeta(\sigma+i\tau t)\right|^2 \textrm{d}\tau + \mathcal{G}(\sigma,t)  \notag \\
&+ O\left(e^{-\pi t^{\delta_{14}}}\right) = 0, \quad 0<\sigma<1, \ \delta_1>0, \ \delta_4>0, \ \delta_{14} = \min{(\delta_1, \delta_4)}, \quad t\to\infty,
\end{align}
where $\Gamma(z)$, $z\in\mathbb{C}$, denotes the gamma function, $\mathcal{G}(\sigma, t)$ is defined by equations \eqref{2.4} and \eqref{2.5b}, and the principal value integral is defined with respect to $\tau=1$.
%\begin{equation} \label{3.3}
%\oint_{-\frac{t^{\delta_1}}{t}}^{1+\frac{t^{\delta_4}}{t}}  \textrm{d}\tau = \lim_{\varepsilon\to 0} \left( \int_{-\frac{t^{\delta_1}}{t}}^{1-\varepsilon}  \textrm{d}\tau + \int_{1+\varepsilon}^{1+\frac{t^{\delta_4}}{t}}  \textrm{d}\tau  \right).
%\end{equation}
\end{lemma}

\textbf{Proof} Starting with Stirling's formula, the following formulae are derived
in the appendix of \cite{FL}:
\addtocounter{equation}{1}
\begin{align} \label{2.19a}
\Gamma(\sigma+i\xi) = \sqrt{2\pi}\xi^{\sigma-\frac{1}{2}}e^{-\frac{\pi \xi}{2}}e^{-\frac{i\pi}{4}}e^{-i\xi}\xi^{i\xi}e^{\frac{i\pi \sigma}{2}} \left[ 1 + O\left(\frac{1}{\xi}\right) \right], \quad \xi \to\infty, \tag{\theequation a}
\end{align}
and
\begin{align} \label{2.19b}
\Gamma(\sigma-i\xi) = \sqrt{2\pi}\xi^{\sigma-\frac{1}{2}}e^{-\frac{\pi \xi}{2}}e^{\frac{i\pi}{4}}e^{i\xi}\xi^{-i\xi}
e^{-\frac{i\pi \sigma}{2}} \left[ 1 + O\left(\frac{1}{\xi}\right) \right], \quad \xi \to\infty. \tag{\theequation b}
\end{align}
These formulae imply that the kernel $K$ defined in \eqref{2.30} satisfies the following estimate:
\begin{align*}
K(\sigma, t, \tau) &=   O\left((|1-\tau|t)^{-\frac{1}{2}} e^{-\frac{\pi}{2}|1-\tau|t} ~ \frac{(|\tau|t)^{\sigma-\frac{1}{2}} e^{-\frac{\pi}{2}|\tau|t}}{t^{\sigma-\frac{1}{2}} e^{-\frac{\pi}{2}t}} \right) \\
&= O\left((|1-\tau|t)^{-\frac{1}{2}} |\tau|^{\sigma-\frac{1}{2}} e^{-\frac{\pi}{2}t (|1-\tau|+|\tau| -1)} \right),
\end{align*}
provided that as $t\to\infty$,
\begin{equation} \label{3.4}
|\tau t|\to\infty, \quad |1-\tau|t\to\infty.
\end{equation}

In order to ensure the validity of \eqref{3.4}, the boundaries $\tau=0$ and $\tau=1$ must be analysed carefully. In this connection, we decompose the infinite line as the union of the following three subintervals:
\begin{equation*}
-\infty<\tau\le- t^{\delta_1-1} , \qquad -t^{\delta_1-1}\le\tau\le 1+ t^{\delta_4-1}, \qquad 1+ t^{\delta_4-1}<\tau< \infty,
\end{equation*}
with $\delta_1>0$, $\delta_4>0$.

For $\tau$ in the first interval we find
\begin{equation*}
e^{-\frac{\pi}{2}t (|1-\tau|+|\tau| -1)} = e^{-\frac{\pi}{2}t (1-\tau-\tau-1)} = e^{\pi\tau t } .
%\le e^{-\pi t^{\delta_1}}.
\end{equation*}
Hence,
\begin{equation}
\Bigg| t\int_{-\infty }^{-\frac{t^{\delta _{1}}}{t}}K(\sigma ,t,\tau )|\zeta (\sigma +it\tau )|^{2}d\tau \Bigg| \leq
t\int_{-\infty }^{-\frac{t^{\delta _{1}}}{t}}\frac{1}{\sqrt{t}}|\tau |^{\sigma -\frac{1}{2}}e^{\pi\tau t}
|\zeta (\sigma +it\tau )|^{2}d\tau  \nonumber
\end{equation}
\begin{equation}
=t^{-\sigma }\int_{t^{\delta _{1}}}^{\infty}x^{\sigma -\frac{1}{2}}e^{-\pi x}|\zeta (\sigma +ix)|^{2}dx
\leq e^{-\pi t^{\delta _{1}}}.  \nonumber
\end{equation}

Similar considerations apply for $\tau$ in the third interval, where
\begin{equation*}
e^{-\frac{\pi}{2}t (|1-\tau|+|\tau| -1)} = e^{-\frac{\pi}{2}t (\tau- 1 + \tau-1)} = e^{-\pi t (\tau -1)}.
%\le e^{-\pi t^{\delta_4}}.
\end{equation*}
Thus, for large $t$ the kernel $K$ is exponentially small for $\tau$ outside the interval $-t^{\delta_1 -1}\le\tau\le 1 + t^{\delta_4 -1}$, and hence equation \eqref{2.23} becomes equation \eqref{3.1}.

\textbf{QED}

\vphantom{a}

 Theorem \ref{t2.1} implies the following result.

\begin{corollary} \label{t3.1}
The Riemann zeta function $\zeta(s)$ satisfies the integral equation
\begin{align} \label{3.6}
&\frac{t}{\pi} \oint_{-t^{\delta_1-1}}^{1+t^{\delta_4-1}} \Re {\left\{ \frac{\Gamma(it-i\tau t)}{\Gamma(\sigma + it)}\Gamma(\sigma + i\tau t) \right\}} \left| \zeta(\sigma+i\tau t)\right|^2 \textrm{d}\tau \notag \\
&+ \begin{cases}
      \zeta(2\sigma) + 2\Gamma(2\sigma - 1) \zeta(2\sigma - 1) \sin{(\pi\sigma) t^{1-2\sigma}} \left(  1+O\left( \frac{1}{t}  \right) \right), & 0<\sigma<1, \quad \sigma \ne \frac{1}{2}, \notag\\
      \ln{t} + 2\gamma - \ln{2\pi}, & \sigma = \frac{1}{2},
   \end{cases}\\
   &+ O \left(  e^{-\pi t^{\delta_{14}}} \right) + O \left( \frac{1}{t^2}\right) = 0, \quad  \delta_1>0, \ \delta_4>0, \ \delta_{14} = \min{(\delta_1, \delta_4)}, \quad t\to\infty,
\end{align}
where $\Gamma(z)$, $z\in\mathbb{C}$, denotes the gamma function, the principal value integral is defined with respect to $\tau=1$.
\end{corollary}

\textbf{Proof}
Employing (2.33) with $\xi=t$ in the definition of $\mathcal{G}(\sigma,t)$ we find
\begin{equation}  \label{2.20}
\mathcal{G}(\sigma,t) = 2 \Gamma(2\sigma - 1)\zeta(2\sigma - 1)\sin{(\pi\sigma)}t^{1-2\sigma}\left[1 + O\left(\frac{1}{t}\right) \right] +\zeta(2\sigma) + O\left(\frac{1}{t^2}\right),  \nonumber
\end{equation}
\begin{equation}
 0<\sigma< 1, ~ t\to\infty.
\end{equation}
Furthermore, the estimate
\begin{equation}  \label{2.21}
\Psi\left(\frac{1}{2} + it\right) = \ln{t} + \frac{i\pi}{2}+ O\left(\frac{1}{t^2}\right), \quad t\to\infty,
\end{equation}
implies
\begin{equation}  \label{2.22}
\mathcal{G}\left(\frac{1}{2}, t\right) = \ln{t} + 2\gamma - \ln{2\pi} + O\left(\frac{1}{t^2}\right), \quad t\to\infty.
\end{equation}
Replacing in equation \eqref{3.1} the function $\mathcal{G}(\sigma,t)$ by equations \eqref{2.20} and \eqref{2.22}, equation \eqref{3.1} becomes equation \eqref{3.6}.

\textbf{QED}

\section{A Relation Between $S_R$ and $S_M$ } \label{sec6}

In what follows we prove equations \eqref{1.21} and \eqref{1.22} relating the double exponential sums $S_R$ and $S_M$ .

\begin{lemma}
Define the functions $f(u,v)$ and $g(u,v)$ by
\begin{equation}  \label{sevenone}
f(u,v)=\sum_{m_{1}=1}^{N}\sum_{m_{2}=1}^{N}\frac{1}{m_{2}^{u}}\frac{1}{(m_{1}+m_{2})^{v}},
\end{equation}
\begin{equation}  \label{seventwo}
g(u,v)=\sum_{m=1}^{N}\sum_{n=N+1}^{N+m}\frac{1}{m^{u}n^{v}},
\end{equation}
where $N$ is an arbitrary finite positive integer and $u\in\mathbb{C}$, $v\in\mathbb{C}$. These functions satisfy the identity
\begin{equation}  \label{seventhree}
f(u,v)+f(v,u)+\sum_{m=1}^{N}\frac{1}{m^{u+v}}=
\bigl( \sum_{m=1}^{N}\frac{1}{m^{u}} \bigr)\bigl( \sum_{n=1}^{N}\frac{1}{n^{v}} \bigr) +g(u,v)+g(v,u).
\end{equation}
\end{lemma}

\textbf{Proof}
Letting $m_{2}=m$, $m_{1}+m_{2}=n$ in $f(u,v)$ and in $f(v,u)$, and then exchanging $m$ and $n$ in the expression of
$f(v,u)$, we find the following:
\begin{align*}
f(u,v)&+f(v,u)=\Bigl( \sum_{m=1}^{N}\sum_{n=m+1}^{m+N}+\sum_{n=1}^{N}\sum_{m=n+1}^{N+n} \Bigr)\frac{1}{m^{u}n^{v}} \\
&=\Bigl( \sum_{m=1}^{N}\sum_{n=m+1}^{N}+\sum_{m=1}^{N}\sum_{n=N+1}^{N+m}
+\sum_{n=1}^{N}\sum_{m=n+1}^{N}+\sum_{n=1}^{N}\sum_{m=N+1}^{N+n}    \Bigr)\frac{1}{m^{u}n^{v}}.  \nonumber
\end{align*}
The second sum above equals $g(u,v)$, and by exchanging $m$ and $n$ in the last sum it follows that the latter sum equals $g(v,u)$.
Thus, the above identity becomes
\begin{equation} \label{sevenfour}
f(u,v)+f(v,u)=\Bigl( \sum_{m=1}^{N}\sum_{n=m+1}^{N}+\sum_{n=1}^{N}\sum_{m=n+1}^{N} \Bigr)\frac{1}{m^{u}n^{v}}
+g(u,v)+g(v,u).
\end{equation}
But
\begin{equation}  \label{sevenfive}
\sum_{n=1}^{N}\sum_{m=n+1}^{N}\frac{1}{m^{u}n^{v}}=\sum_{m=1}^{N}\sum_{n=1}^{m-1}\frac{1}{m^{u}n^{v}}.
\end{equation}

Using the identity (\ref{sevenfive}) in (\ref{sevenfour}), adding to both sides of (\ref{sevenfour}) the term
\begin{equation}
\sum_{m=1}^{N}\frac{1}{m^{u}m^{v}},  \nonumber
\end{equation}
and noting that
$$\Bigl( \sum_{m=1}^{N}\sum_{n=m+1}^{N}+\sum_{m=1}^{N}\sum_{n=1}^{m-1}\Bigr)  \frac{1}{m^{u}n^{v}}
+\sum_{m=1}^{N}\frac{1}{m^{u}{m}^{v}}=\sum_{m=1}^{N}\sum_{n=1}^{N}\frac{1}{m^{u}n^{v}}, $$
equation (\ref{sevenfour}) becomes (\ref{seventhree}).

\textbf{QED}
\begin{corollary} \label{c7.1}The following identity is valid:
\begin{align}\label{sevensix}
2\Re&\left\{ 
\sum_{m_{1}=1}^{[t]}\sum_{m_{2}=1}^{[t]}\frac{1}{m_{2}^{\bar{s}}(m_{1}+m_{2})^{s }} \right\} -
\left( \sum_{m=1}^{[t]}\frac{1}{m^{s}} \right) \left( \sum_{m=1}^{[t]}\frac{1}{m^{\bar{s} }} \right)\nonumber \\ 
&=-\sum_{m=1}^{[t]}\frac{1}{m^{2\sigma }}
+2\Re \left\{ \sum_{m=1}^{[t]}\sum_{n=[t]+1}^{[t]+m}\frac{1}{m^{\bar{s} }n^{s}} \right\}, 
\qquad s=\sigma + i t.
\end{align}
\end{corollary}

\textbf{Proof} Letting $u=\bar{s}$, $v=s $, $N=[t]$, equation \eqref{seventhree} becomes\eqref{sevensix}.

\textbf{QED}

%It is shown in \cite{F} that the last two terms in the rhs of  (1.13) are ``small" as $t\rightarrow\infty$. Thus, in
%what follows we concentrate on the comparison on the lhs of (1.13) with the sum in (1.9).

\begin{lemma} \label{l7.3}
The following identity is valid: 
\begin{equation}   \label{sevenfifteen}
\sum_{m_{1}=1}^{[t]}\sum_{m_{2}=1}^{[t]}=\sum_{(m_{1},m_{2})\in M}+
\sum_{m_{1}=1}^{\left[\frac{t}{t^{1-\delta _{3}}-1}\right]-1}
\sum_{m_{2}=\left[(t^{1-\delta _{3}}-1)m_{1}\right]+1}^{[t]}+
\sum_{m_{1}=\left[t^{1-\delta _{2}}\right]}^{[t]}\sum_{m_{2}=1}^{\left[\frac{m_{1}}{t^{1-\delta _{2}}-1}\right]-1},
\end{equation}
where the set $M$ is defined in \eqref{1.14}.
\end{lemma}

\textbf{Proof} 
\begin{equation}
\sum_{m_{1}=1}^{[t]}\sum_{m_{2}=1}^{[t]}=
\Biggl(  
\sum_{m_{1}=1}^{\left[\frac{t}{t^{1-\delta _{3}}-1}\right]-1}
+\sum_{m_{1}=\left[\frac{t}{t^{1-\delta _{3}}-1}\right]}^{\left[t^{1-\delta _{2}}\right]-1}
+\sum_{m_{1}=\left[t^{1-\delta _{2}}\right]}^{[t]}
\Biggr)\sum_{m_{2}=1}^{[t]}  \nonumber
\end{equation}
\begin{equation}   \label{sevensixteen}
=\sum_{m_{1}=1}^{\left[\frac{t}{t^{1-\delta _{3}}-1}\right]-1}\sum_{m_{2}=1}^{[t]}
+\sum_{m_{1}=\left[\frac{t}{t^{1-\delta _{3}}-1}\right]}^{\left[t^{1-\delta _{2}}\right]-1}\sum_{m_{2}=1}^{[t]}
+\sum_{m_{1}=\left[t^{1-\delta _{2}}\right]}^{[t]}\sum_{m_{2}=1}^{[t]}.
\end{equation}
We subdivide the sum over $m_{2}$ occurring in the first and third double sums in the second equality of (\ref{sevensixteen}) as follows:

\begin{equation}
\sum_{m_{2}=1}^{[t]}=\sum_{m_{2}=1}^{\left[(t^{1-\delta _{3}}-1)m_{1}\right]}+
\sum_{m_{2}=\left[(t^{1-\delta _{3}}-1)m_{1}\right]+1}^{[t]}, \nonumber
\end{equation}
and
\begin{equation}
\sum_{m_{2}=1}^{[t]}= \sum_{m_{2}=1}^{\left[\frac{m_{1}}{t^{1-\delta _{2}}-1}\right]-1}
+\sum_{m_{2}=\left[\frac{m_{1}}{t^{1-\delta _{2}}-1}\right]}^{[t]}.  \nonumber
\end{equation}
Substituting the above expressions in (\ref{sevensixteen}) we find
\begin{equation}
\sum_{m_{1}=1}^{[t]}\sum_{m_{2}=1}^{[t]}=
\sum_{m_{1}=1}^{\left[\frac{t}{t^{1-\delta _{3}}-1}\right]-1}
\Biggl(
\sum_{m_{2}=1}^{\left[(t^{1-\delta _{3}}-1)m_{1}\right]}+
\sum_{m_{2}=\left[(t^{1-\delta _{3}}-1)m_{1}\right]+1}^{[t]}
\Biggr)   \nonumber
\end{equation}
\begin{equation}  \label{sevenseventeen}
+\sum_{m_{1}=\left[\frac{t}{t^{1-\delta _{3}}-1}\right]}^{\left[t^{1-\delta _{2}}\right]-1}\sum_{m_{2}=1}^{[t]}+
\sum_{m_{1}=\left[t^{1-\delta _{2}}\right]}^{[t]}
\Biggl(
\sum_{m_{2}=1}^{\left[\frac{m_{1}}{t^{1-\delta _{2}}-1}\right]-1}+
\sum_{m_{2}=\left[\frac{m_{1}}{t^{1-\delta _{2}}-1}\right]}^{[t]}
\Biggr).
\end{equation}
The sum of the first, third, and fifth double sums in (\ref{sevenseventeen}) equals the first term of the rhs of 
(\ref{sevenfifteen}), whereas the second and fourth double sums in (\ref{sevenseventeen}) are the second and third 
terms of the rhs of (\ref{sevenfifteen}).

\textbf{QED}

The sums appearing in the rhs of \eqref{sevensix} can be estimated using the results of the lemma below.

\begin{lemma}\label{lemma6.3n}
\begin{align}\label{seventen-new}
2\Re&\left\{ 
\sum_{m_{1}=1}^{[t]}\sum_{m_{2}=1}^{[t]}\frac{1}{m_{2}^{\bar{s}}(m_{1}+m_{2})^{s }} \right\} -
\left| \sum_{m=1}^{[t]}\frac{1}{m^{s}}  \right|^2\nonumber \\
&=-\begin{cases}\ln t+O(1), \quad &\sigma=\frac{1}{2}, \vspace*{2mm} \\
\dfrac{t^{1-2\sigma }}{1-2\sigma }+O(1), \quad &0<\sigma <1, \quad \sigma \neq\frac{1}{2}, \end{cases}\notag\\
&+\begin{cases}
O \left(  t^{\frac{1}{2}-\frac{5}{3}\sigma} \ln{t} \right), & 0<\sigma\leq \frac{1}{2},\vspace*{2mm} \\
O \left(  t^{\frac{1}{3}-\frac{4}{3}\sigma} \ln{t}  \right), & \frac{1}{2} < \sigma< 1,
\end{cases} \qquad t\to \infty.
\end{align}
\end{lemma}

\textbf{Proof} 
 The first term of the rhs of \eqref{sevensix} yields $$\sum_{m=1}^{[t]}\frac{1}{m^{2\sigma }}=\int_0^{[t]}\frac{1}{x^{2\sigma}}dx + O(1),$$
which gives rise to the first term of the rhs of \eqref{seventen-new}. The second term of the rhs of \eqref{sevensix} is estimated in Lemma 3.2 of \cite{KF}, and then \eqref{seventen-new} follows.

\textbf{QED}

\section{The Rigorous Estimation of $I_1$ and $I_2$}

\begin{lemma} \label{lemma31}
Let $I_1(\sigma, t, \delta_1)$ be defined by
\begin{equation} \label{eq31}
I_1(\sigma, t, \delta_1) = \frac{t}{\pi} \int_{-t^{\delta_1 -1}}^{\frac{1}{t}} \Re{ \left\{\frac{\Gamma(it-it\tau)}{\Gamma(\sigma+it)} \Gamma(\sigma+it\tau) \right\} } |\zeta(\sigma+it\tau)|^2 d\tau,
\end{equation}
\begin{equation*}
0<\sigma<1, \quad t>0, \quad 0< \delta_1 <1.
\end{equation*}
Then, for sufficiently small $\delta_1$,
\begin{equation} \label{eq32}
I_1(\sigma, t, \delta_1) = \times\begin{cases}
O \left(  \frac{t^{\left(2-\frac{4}{3}\sigma\right)\delta_1}}{t^\sigma} \right), & 0\leq\sigma\leq \frac{1}{2},\\
O \left(  \frac{t^{\left(\frac{5}{3}-\frac{2}{3}\sigma\right)\delta_1}} {t^\sigma}  \right), & \frac{1}{2} < \sigma< 1,
   \end{cases}, \quad t \to \infty.
\end{equation}
\end{lemma}

\noindent \textbf{Proof} In the interval of integration, we have
\begin{equation} \label{eq33}
-t^{\delta_1}\leq t \tau \leq 1.
\end{equation}
Thus,
\begin{equation*}
t-1 \leq t -\tau t \leq t + t^{\delta_1}.
\end{equation*}
Hence, $t-\tau t \to \infty$ as $t \to \infty$, and therefore we can use the asymptotic formulas (2.33a) to compute both $\Gamma(\sigma+it)$ and $\Gamma(it - it\tau)$:
\begin{subequations} \label{eq34}
\begin{equation} \label{eq34a}
\Gamma(\sigma+it) = \sqrt{2\pi} t^{\sigma-\frac{1}{2}}e^{-\frac{\pi t}{2}}e^{-\frac{i \pi}{4}} e^{\frac{i \pi \sigma}{2}}e^{-it} t^{it} \left[  1 + O \left(  \frac{1}{t} \right) \right], \quad t \to \infty,
\end{equation}
\begin{multline} \label{eq34b}
\Gamma(it - it\tau) = \sqrt{2\pi} t^{-\frac{1}{2}} (1-\tau)^{-\frac{1}{2}} e^{-\frac{\pi t}{2}(1-\tau)}e^{-\frac{i \pi}{4}} e^{-it(1-\tau)} t^{it(1-\tau)} (1-\tau)^{it(1-\tau)} \\ \times \left[  1 + O \left(  \frac{1}{t-t\tau} \right) \right], \quad t \to \infty.
\end{multline}
\end{subequations}
Equations \eqref{eq34} together with the inequality
\begin{equation*}
\frac{1}{t-\tau t} \leq \frac{1}{t-1},
\end{equation*}
imply that
\begin{equation} \label{eq35}
\frac{\Gamma(it-it\tau)}{\Gamma(\sigma+it)} = t^{-\sigma} (1-\tau)^{-\frac{1}{2}} e^{\frac{\pi t \tau}{2}} e^{i\varphi(t,\tau)} \left[  1 + O \left(  \frac{1}{t} \right) \right], \quad t \to \infty,
\end{equation}
where $\varphi$ is a real function.

Equation \eqref{eq33} implies that $\tau t$ is either O(1), or $\tau t = O (t^{\Delta}), \ \Delta<\delta_1$. 

If $\tau t = O(1)$, then
\begin{equation*}
\Gamma(\sigma+i\tau t) = O(1), \quad t \to \infty.
\end{equation*}
If $\tau t\to -\infty$, then using (2.33b) with $\xi=\tau t$, we find
\begin{equation*}
e^{\frac{\pi\tau t}{2}}\Gamma(\sigma+i\tau t) = O \left(  \left| \tau t\right|^{\sigma-\frac{1}{2}}e^{\pi\tau t} \right) , \quad \tau t \to -\infty.
\end{equation*}
The lhs of the above equation is bounded and exponentially decaying  as $\tau t \to -\infty$.
Thus, the last two estimates yield 
\begin{equation*}
e^{\frac{\pi\tau t}{2}} \Gamma(\sigma+i\tau t)= 
     O\left( 1 \right),   \qquad 0<\sigma <1.
\end{equation*}

Combining the above equations with \eqref{eq35}, we find
\begin{align}\label{eqnew}
\Re{ \left\{\frac{\Gamma(it-it\tau)}{\Gamma(\sigma+it)} \Gamma(\sigma+it\tau) \right\} } = &  O \left(  t^{-\sigma} \right) \left(  1 - \frac{O(1)}{t} \right)^{-\frac{1}{2}}  \left[  1 + O \left( \frac{1}{t}  \right) \right], \notag\\
& 0<\sigma<1, \quad t \to \infty.
\end{align}

Theorem 5.12 in \cite{T} states that
\begin{equation*}
\zeta(1/2 + it) = O \left(  t^{\frac{1}{6}} \ln t \right), \quad t \to \infty,
\end{equation*}
which is improved by Theorem 5.18 in \cite{T} which states the slightly stronger result
\begin{equation*}
\zeta(1/2 + it) = O \left(  t^{\frac{1}{6}}  \right), \quad t \to \infty.
\end{equation*}
The Phragm\'en-Lindel{\"o}f convexity principle (PL)  (known also as Lindel{\"o}f's theorem) implies
\begin{equation*} \label{est0}
\zeta(\sigma + it) = \begin{cases}
O \left(  t^{\frac{1}{2}-\frac{2}{3}\sigma}  \right), & 0\leq\sigma\leq \frac{1}{2},\\
O \left(  t^{\frac{1}{3}-\frac{1}{3}\sigma}   \right), & \frac{1}{2} < \sigma< 1.
\end{cases}
\end{equation*}

The exponents $\ \frac{1}{2}-\frac{2}{3}\sigma \ $ and $\ \frac{1}{3}-\frac{1}{3}\sigma \ $ have been improved only slightly with the best current result due to Bourgain \cite{B}.

The above estimate yields
\begin{equation} \label{eqnew2}
\zeta(\sigma + it^{\Delta}) =\begin{cases}
O \left(  t^{\left(\frac{1}{2}-\frac{2}{3}\sigma\right)\Delta}  \right), & 0\leq\sigma\leq \frac{1}{2},\\
O \left(  t^{\left(\frac{1}{3}-\frac{1}{3}\sigma\right)\Delta}   \right), & \frac{1}{2} < \sigma< 1,
\end{cases}\quad t \to \infty,
\end{equation}
where $\Delta$ is a positive constant. Combining the estimates \eqref{eqnew} and \eqref{eqnew2} we find
\begin{align*}
\Re{ \left\{\frac{\Gamma(it-it\tau)}{\Gamma(\sigma+it)} \Gamma(\sigma+it\tau) \right\} } |\zeta(\sigma+it\tau)|^2 &= O \left(t^{- \sigma}\right)  \left[  1 + O \left( \frac{1}{t}  \right) \right]\\
&\times\begin{cases}
O \left(  t^{\left(1-\frac{4}{3}\sigma\right)\delta_1}  \right), & 0\leq\sigma\leq \frac{1}{2},\\
O \left(  t^{\left(\frac{2}{3}-\frac{2}{3}\sigma\right)\delta_1}   \right), & \frac{1}{2} < \sigma< 1,
   \end{cases}, \quad t \to \infty.
\end{align*}
The above equation together with the mean value theorem, imply
\begin{align*}
I_1 = \frac{t}{\pi} \left( \frac{1}{t} + \frac{t^{\delta_1}}{t} \right) O \left(t^{- \sigma}\right) \times\begin{cases}
O \left(  t^{\left(1-\frac{4}{3}\sigma\right)\delta_1}  \right), & 0\leq\sigma\leq \frac{1}{2},\\
O \left(  t^{\left(\frac{2}{3}-\frac{2}{3}\sigma\right)\delta_1}   \right), & \frac{1}{2} < \sigma< 1,
   \end{cases}, \quad t \to \infty.
\end{align*}
which yields \eqref{eq32}.
\textbf{QED}

\begin{lemma} \label{lemma32}
Let $I_2(\sigma, t, \delta_2)$ be defined by
\begin{equation} \label{eq36}
I_2(\sigma, t, \delta_2) = \frac{t}{\pi} \int_{\frac{1}{t}}^{t^{\delta_2 -1}} \Re{ \left\{\frac{\Gamma(it-it\tau)}{\Gamma(\sigma+it)} \Gamma(\sigma+it\tau) \right\} } |\zeta(\sigma+it\tau)|^2 d\tau,
\end{equation}
\begin{equation*}
0<\sigma<1, \quad t>0, \quad 0< \delta_2 <1.
\end{equation*}
Then,
\begin{equation} \label{eq37}
I_2\left(\frac{1}{2}, t,\delta_2\right) = O\left( t^{-\frac{1}{2}+ \delta_2} \ln{t} \right), \quad \sigma = \frac{1}{2},\quad t\to\infty,
\end{equation}
and
\begin{equation} \label{eq38}
I_2(\sigma, t,\delta_2) =  \begin{cases}
     O\left( t^{-\sigma + \left(\sigma + \frac{1}{2} \right)\delta_2}\zeta(2\sigma)\right), &  \frac{1}{2}<\sigma< 1, \quad t\to\infty,\\
     O\left( t^{-\sigma + 2 (1 -\sigma )\delta_2} \zeta(2-2\sigma) \right), &  0<\sigma<\frac{1}{2}, \quad t\to\infty.
   \end{cases}
\end{equation}
\end{lemma}
\noindent \textbf{Proof} In the interval of integration, we have
\begin{equation*}
1\le t\tau\le t^{\delta_2}.
\end{equation*}
Thus,
\begin{equation*}
t - t^{\delta_2} \leq t -t \tau \leq t -1.
\end{equation*}
Hence, $t-t\tau\to\infty$ as $\tau\to\infty$, and therefore we can use (2.33a) to compute $|\Gamma(it - it\tau)$. Thus, taking into consideration the inequality
\begin{equation*}
\frac{1}{t-t\tau} \leq \frac{1}{t - t^{\delta_2}} = \frac{1}{t(1-t^{\delta_2-1})} = \frac{1}{t} + O( t^{\delta_2 -2}), \quad t \to \infty,
\end{equation*}
we find that equation \eqref{eq35} is still valid. Inserting the expression \eqref{eq35} in the definition of $I_2$ and using the change of variables $\rho=t\tau$, we find
\begin{multline} \label{eq39}
I_2 = \frac{1}{\pi t^{\sigma}} \int_{1}^{t^{\delta_2}} \left(  1 - \frac{\rho}{t} \right)^{-\frac{1}{2}} \Re{\left\{ \left( e^{\frac{\pi\rho}{2}} \Gamma(\sigma+i\rho) \right) e^{i\varphi(t,\rho)} \right\}} |\zeta(\sigma+i\rho)|^2 d\rho \\ \times \left[ 1 + O \left(  \frac{1}{t}  \right) \right], \quad t \to \infty.
\end{multline}
We will estimate the integral \eqref{eq39} by employing the first mean value theorem for integrals: since $\left| \zeta(\sigma+i\rho)\right|^2$ does not change sign and it is integrable for $\rho\in[1,t^{\delta_2}]$, it follows that there exists a $c(t)$,
\begin{equation*}
1<c(t)<t^{\delta_2},
\end{equation*}
such that
\begin{multline} \label{eq310}
I_2= \frac{1}{\pi t^{\sigma}} \left(  1 - \frac{c(t)}{t} \right)^{-\frac{1}{2}} \Re{\left\{ \left( e^{\frac{\pi c(t)}{2}} \Gamma(\sigma+ic(t)) \right) e^{i\varphi(t,c(t))} \right\}} \int_{1}^{t^{\delta_2}} |\zeta(\sigma+i\rho)|^2 d\rho \\ \times \left[ 1 + O \left(  \frac{1}{t}  \right) \right], \quad t \to \infty.
\end{multline}
If $c(t)=O(1)$, then $e^{\frac{\pi c(t)}{2}} \Gamma(\sigma+ic(t))=O(1)$. If $c(t)\to\infty$, as $t\to\infty$, then equation (2.33a) with $\xi=c(t)$ yields
\begin{equation*}
e^{\frac{\pi c(t)}{2}} \Gamma(\sigma+ic(t)) = O \left( c(t)^{ \left(  \sigma -\frac{1}{2} \right)}  \right) ,\quad t \to \infty.
\end{equation*}
Thus, since $c(t)<t^{\delta_2}$, we find
\begin{equation} \label{eq311}
e^{\frac{\pi c(t)}{2}} \Gamma(\sigma+ic(t))= \begin{cases}
     O\left( t^{\delta_2 \left( \sigma - \frac{1}{2} \right)} \right), &  \frac{1}{2}<\sigma< 1,\\
     O\left( 1 \right), &  0<\sigma \leq \frac{1}{2}.
   \end{cases}
\end{equation}
We recall Atkinson's asymptotic formula (Theorem 7.4 of \cite{T}):
\begin{align} \label{eq312}
\int_{1}^{t^{\delta_2}} \left| \zeta\left(\frac{1}{2}+i\rho\right)\right|^2 \textrm{d}\rho =& t^{\delta_2} \ln{t^{\delta_2}} + (2\gamma - 1 - \ln{2\pi}) t^{\delta_2}  \nonumber \\ &+ O\left(  t^{\delta_2 \left( \frac{1}{2} + \varepsilon \right)}\right), 
\qquad \varepsilon >0, \quad t\to\infty.
\end{align}
Similarly,  for $\sigma\ne\frac{1}{2}$ we obtain the following estimate:
\begin{equation} \label{eq313}
\int_{1}^{t^{\delta_2}} \left| \zeta(\sigma+i\rho)\right|^2 \textrm{d}\rho = \begin{cases}
    O\left( t^{\delta_2}\zeta(2\sigma)\right), & \frac{1}{2}<\sigma\le 1, \quad  t\to\infty,  \\
    O\left( t^{2(1 -\sigma)\delta_2} \zeta(2 - 2\sigma) \right), &  0<\sigma<\frac{1}{2}, \quad t\to\infty.
   \end{cases}
\end{equation}
The above estimates are obtained using the results of the Theorem 7.2 of \cite{T}. Indeed using these results for $\sigma=1/2$ one obtains \eqref{eq312}, and for 
$1/2<\sigma<1$  one obtains the first of \eqref{eq313}. It is straightforward to make the analogous analysis contained in the proof of Theorem 7.2 of \cite{T} for 
$0<\sigma<1/2$, which in turn yields the second of \eqref{eq313}.

Using in equation \eqref{eq310} the expressions obtained from equations \eqref{eq311}-\eqref{eq313}, we find equation \eqref{eq38}. \textbf{QED}

Combining Lemma 2.2 and Lemmas \ref{lemma31} and \ref{lemma32} we obtain the following result.

\begin{theorem}\label{theorem31}
Let $0<\sigma<1$, $t>0$, and let $\delta_1$ be a sufficiently small positive constant, whereas the constants $\delta_j$ satisfy $0<\delta_j<1, \ j=2,3,4$. Let the integrals $I_3(\sigma, t, \delta_2, \delta_3)$ and $I_4(\sigma, t, \delta_3, \delta_4)$ be defined by (1.8) with $j=3$ and $j=4$. The Riemann zeta function $\zeta(s)$, $s=\sigma+it$ satisfies the following equation:
\begin{multline} \label{eq314}
I_{3}(\sigma ,t,\delta _{2},\delta _{3}) + I_4(\sigma, t,\delta_3,\delta_4)  +     \ln{t} + 2\gamma - \ln{2\pi} + O \left(  t^{\delta_2-\frac{1}{2}} \ln{t} \right)+ O \left(  t^{\frac{4}{3}\delta_1-\frac{1}{2}} \right)\\
   + O \left(  e^{-\pi t^{\delta_{14}}} \right) = 0, \quad \delta_{14}=\min{(\delta_1, \delta_4)}, \quad \sigma = \frac{1}{2}, \quad t\to\infty,
\end{multline}
as well as
\begin{multline} \label{eq315}
I_{3}(\sigma ,t,\delta _{2},\delta _{3}) + I_4(\sigma, t,\delta_3,\delta_4)+ \zeta(2\sigma) + 2\Gamma(2\sigma - 1) \zeta(2\sigma - 1) \sin{(\pi\sigma) t^{1-2\sigma}} \left(  1+O\left( \frac{1}{t}  \right) \right) \\
+ \begin{cases}
     O\left( t^{-\sigma + \left(\frac{3}{2} -\sigma  \right)\delta_2} \ln{t} \right)+O \left(  t^{-\sigma+\frac{4}{3}\delta_1} \right), & 0<\sigma < \frac{1}{2} \\
     O\left( t^{-\sigma + \left(\sigma + \frac{1}{2} \right)\delta_2}\zeta(2\sigma)\right)+ O\left( t^{-\sigma + \left(\frac{5}{6} +\sigma  \right)\delta_1}\right), &  \frac{1}{2} < \sigma <1
   \end{cases}\\
+ O \left(  e^{-\pi t^{\delta_{14}}} \right) = 0, \quad t\to\infty.
\end{multline}
\end{theorem}
\noindent \textbf{Proof} Decomposing the interval of integration of the lhs of equation (2.35) into the four subintervals defined in \eqref{1.6}, we find that the lhs of equation (2.35) equals  the sum of the four integrals $\{I_{j}\}_{1}^{4}$ defined by equations (1.8). Replacing the integrals $I_1$ and $I_2$ via equations \eqref{eq31} and \eqref{eq36} respectively, equation (2.35) yields equations \eqref{eq314} and \eqref{eq315}. \textbf{QED}

\section{The Asymptotics of $\tilde{I}_3$}

Equation \eqref{1.10} expresses $\tilde{I}_3$ in terms of the integral $J_3$, which we analyse below.

\begin{proposition} \label{prop41}
Let $J_3$ be defined by
\begin{multline} \label{eq41}
J_3(\sigma,t,\delta_2, \delta_3, \lambda) = \frac{t}{\pi} \int_{t^{\delta_2 -1}}^{1-t^{\delta_3 -1}} \frac{\Gamma(it-it\tau)}{\Gamma(\sigma+it)} \Gamma(\sigma+it\tau) \lambda^{i\tau t }  d\tau, \quad \lambda = \frac{m_2}{m_1}, \\ 0 < \sigma <1, \quad t>0, \quad 0 < \delta_2 <1, \ 0<\delta_3<1, \quad m_j = 1, 2, \dots, [t], \quad j= 1,2.
\end{multline}
Then,
\begin{multline} \label{eq42}
J_3(\sigma,t,\delta_2, \delta_3, \lambda) = \sqrt{\frac{2t}{\pi}}  e^{-\frac{i\pi}{4}} \tilde{J}_3(\sigma,t,\delta_2, \delta_3, \lambda)   \left[ 1 + O(t^{-\delta_{23}})  \right], \\ \delta_{23} = \min{\left\{ \delta_2, \delta_3  \right\}}, \quad t \to\infty,
\end{multline}
where $\tilde{J}_3$ is defined by
\begin{equation} \label{eq43}
\tilde{J}_3 (\sigma,t,\delta_2, \delta_3, \lambda) = \int_{t^{\delta_2 -1}}^{1-t^{\delta_3 -1}} G(\sigma, \tau) e^{itF(\tau, \lambda)}d\tau,
\end{equation}
with
\begin{equation} \label{eq44}
G(\sigma, \tau) = (1-\tau)^{-\frac{1}{2}}\tau^{\sigma-\frac{1}{2}}, \quad F(\tau, \lambda) = (1-\tau) \ln{(1-\tau)} + \tau \ln{\tau} + \tau \ln{\lambda}.
\end{equation}
\end{proposition}
\noindent \textbf{Proof} In the interval of integration, we have
\begin{equation} \label{eq45}
t^{\delta_2} \leq t\tau \leq t - t^{\delta_3}.
\end{equation}
Thus,
\begin{equation} \label{eq46}
t^{\delta_3} \leq t - t\tau \leq t - t^{\delta_2}.
\end{equation}
Hence, $t\tau\to\infty$ and $t-t\tau\to\infty$ as $t\to\infty$, therefore we can employ the asymptotic formula (2.33a) with $\xi = t$, $t\tau$, $t-\tau$, to compute the ratio of the gamma functions appearing in the rhs of \eqref{eq41}. Expressions for $\Gamma(\sigma+it)$ and the $\Gamma(it-it\tau)$ are given in \eqref{eq34}. Similarly,
\begin{equation} \label{eq47}
\Gamma(\sigma+it\tau ) = \sqrt{2\pi} (t\tau)^{\sigma-\frac{1}{2}}e^{-\frac{\pi \tau t}{2}}e^{-\frac{i \pi}{4}} e^{\frac{i \pi \sigma}{2}}e^{-i\tau t} (\tau t)^{i\tau t} \left[  1 + O \left(  \frac{1}{\tau t} \right) \right], \quad t \to \infty.
\end{equation}

Equations \eqref{eq45} and \eqref{eq46} imply the inequalities
\begin{equation} \label{eq48}
\frac{1}{\tau t} \leq t^{-\delta_2}, \quad \frac{1}{t - t \tau} \leq t^{-\delta_3}.
\end{equation}
Equations \eqref{eq34}, \eqref{eq47} and \eqref{eq48} yield
\begin{multline} \label{eq49}
\frac{\Gamma(it-i t\tau )}{\Gamma(\sigma+it)} \Gamma(\sigma+it\tau) = \sqrt{\frac{2\pi}{t}}  e^{-\frac{i \pi}{4}} G(\sigma,\tau) e^{it[(1-\tau) \ln{(1-\tau)} + \tau \ln{\tau}]} \\ \times \left[ 1 + O(t^{-\delta_{23}})  \right], \quad t\to\infty.
\end{multline}
Substituting equation \eqref{eq49} in the rhs of \eqref{eq41}, and taking into account that the integral of $G$ is positive, we obtain \eqref{eq42}. \textbf{QED}

It is well known that the main contributions to the asymptotic analysis of integrals such as $\tilde{J}_3$ come from possible singularities, from possible stationary points, and from the end points of the interval of integration [AF]. The integral $\tilde{J}_3$ possesses a stationary point at $\tau =1/(1+\lambda )$.  Thus, for this integral there exist two contributions, one from the associated stationary point and one from the end points of the interval of integration. Assuming that the stationary point does not approach the end points, the latter contributions can be computed using integration by parts
%: the contribution to $\tilde{J}_3$ from the upper and lower end points are of the following order respectively:
 together with the following estimates:
\begin{equation} \label{eq410}
\frac{G(\sigma, 1 -t^{\delta_3 -1})}{t} \sim \frac{t^{-\frac{\delta_3 -1}{2}}}{t} = \frac{t^{-\frac{\delta_3}{2}}}{\sqrt{t} }, \quad t \to \infty,
\end{equation}
and
\begin{equation} \label{eq411}
\frac{G(\sigma, t^{\delta_2 -1})}{t} \sim \frac{t^{(\delta_2 - 1) \left(  \sigma - \frac{1}{2} \right)}}{t} = \frac{t^{\delta_2\left(  \sigma - \frac{1}{2} \right)}}{t^{ \sigma + \frac{1}{2}} }, \quad t \to \infty.
\end{equation}

We first compute the contribution from the stationary point, where we also include the error term arising from the contributions of the lower end point.

\begin{proposition} \label{prop42}
Let $\tilde{J}_3$ be defined by \eqref{eq43}. Then,
\begin{equation} \label{eq412}
\tilde{J}_3 = \tilde{J}_3^S -\tilde{J}_3^U,
\end{equation}
where $\tilde{J}_3^S$ and $\tilde{J}_3^U$ are defined as follows:
\begin{equation} \label{eq413}
 \tilde{J}_3^S (\sigma,t,\delta_2, \lambda) = \int_{L(\delta_2)} G(\sigma,\tau) e^{itF(\tau, \lambda)}d\tau,
\end{equation}
and
 \begin{equation} \label{eq414}
 \tilde{J}_3^U (\sigma,t,\delta_3, \lambda) = \int_{1-t^{\delta_3-1}}^{\infty e^{i\varphi}} G(\sigma,\tau) e^{itF(\tau, \lambda)}d\tau, \quad 0<\varphi < \arctan{\frac{\pi}{|\ln{\lambda}|}},
\end{equation}
$\sigma$, $t$, $\delta_2$, $\delta_3$, $\lambda$ are as in \eqref{eq41}, whereas $L(\delta_2)$ denotes the contour in the complex $\tau$-plane, starting at the point $t^{\delta_2 -1}$, going down into the lower half complex plane, up through the point $\tau=1/(1+\lambda)$ and continuing to $\infty e^{i\varphi}$.

$\tilde{J}_3^S$ is given by
\begin{equation} \label{eq415}
\tilde{J}_3^S (\sigma,t,\delta_2, \lambda) =  \sqrt{\frac{2\pi}{t}}  e^{\frac{i \pi}{4}} \frac{\lambda^{it}}{(1+\lambda)^{\sigma+it}} %+ O \left(  \frac{t^{\delta_2 \left(\sigma-\frac{1}{2}\right)}}{t^{\sigma + \frac{1}{2}}}, \lambda \right)
[1+o(1)], \quad t \to \infty, \tag{5.15a}
\end{equation}
where the first term in \eqref{eq415} occurs iff
\begin{equation}\label{eq415b}
\frac{1}{t^{1-\delta_3}-1} < \lambda < t^{1 - \delta_2}-1, \tag{5.15b}
\end{equation}
\addtocounter{equation}{1}
whereas $\tilde{J}_3^U$ is given by
\begin{multline} \label{eq416}
\tilde{J}_3^U (\sigma,t,\delta_3, \lambda) =  \frac{it^{-\frac{\delta_3}{2}}}{\sqrt{t}} t^{i(\delta_3 -1)t^{\delta_3}} (1-t^{\delta_3 -1}) ^ {\sigma - \frac{1}{2}+ i(t-t^{\delta_3})} \frac{\lambda^{i(t-t^{\delta_3})}}{\ln{ (\lambda (t^{1-\delta_3}-1))}} \\ %+ O \left(  \frac{t^{-\delta_3}}{\sqrt{t}}, \lambda \right)
\times [1+o(1)], \quad t \to \infty.
\end{multline}
\end{proposition}
%The notation $O(A, \lambda)$ indicates that the error terms depend on $\lambda$.

\noindent \textbf{Proof} For the function $F(\tau, \lambda)$, $\tau \in \mathbb{C}$, defined by the second of equations \eqref{eq44}, we chose the branch cuts $[-\infty,0] \cup [1,\infty]$. The function $F$ satisfies the equation
\begin{equation} \label{eq417}
\frac{\partial F}{\partial \tau }=-\ln (1-\tau )+\ln \tau +\ln \lambda .
\end{equation}
Thus, the integral $\tilde{J}_3$ possesses a stationary point at $\tau =\tau _{1}$, provided that $1-\tau_1=\tau_1 \lambda$, i.e.,
\begin{equation} \label{eq418}
\tau _{1}=\frac{1}{1+\lambda }.
\end{equation}
This occurs if and only if
\begin{equation*}
\frac{t^{\delta _{2}}}{t}<\tau _{1}<1-\frac{t^{\delta _{3}}}{t},
\end{equation*}
i.e., if and only if $\lambda $ satisfies the inequality \eqref{eq415b}.

We deform the contour of integration to the contour $L(\delta_2)$ along the steepest descent direction, plus the contour from $\infty e^{i\varphi}$ back to the point $1-t^{\delta _{3}-1}$.

We claim that if $\varphi $ is sufficiently small, namely if $\varphi $ satisfies the second of equations in \eqref{eq414}, then the integral $\tilde{J}_{3}^U$ converges. Indeed, employing the change of variables
\begin{equation*}
\tau =\Delta (t)+\rho e^{i\varphi }, \hphantom{3a} \Delta(t) =1-\frac{t^{\delta _{3}}}{t},
\end{equation*}
we find that $F$ becomes
\begin{equation*}
F = (1-\Delta -e^{i\varphi }\rho )\ln (1-\Delta -e^{i\varphi }\rho ) + (\Delta +e^{i\varphi }\rho )[\ln\lambda +\ln (\Delta +e^{i\varphi }\rho )].
\end{equation*}
For fixed $t$ and large $\rho$, we have
\begin{equation}
F\sim \rho e^{i\varphi }[\ln \lambda +\ln (\rho e^{i\varphi })-\ln (-\rho e^{i\varphi })], \hphantom{3a}
\rho \rightarrow \infty .   \nonumber
\end{equation}
Using
\begin{equation*}
\ln{(\rho e^{i\varphi})} - \ln{(-\rho e^{i\varphi})} = i\pi,
\end{equation*}
it follows that
\begin{equation*}
\Im F\sim \rho [\ln\lambda \sin\varphi +\pi\cos\varphi ],  \hphantom{3a} \rho \rightarrow \infty .
\end{equation*}
For the convergence of $\tilde{J}_{3}^U$ we require $\Im F>0$ as $\rho \rightarrow \infty $. Thus, if $\lambda \geq 1$ then we have convergence for all $\varphi \in (0,\pi/2)$, whereas if $\lambda \in (0,1)$ we require that $\varphi $ satisfies the condition displayed in \eqref{eq414}.

In order to compute the contribution from the stationary point $\tau =\tau _{1}$ we employ the well known formula [M]
\begin{multline}  \label{eq419}
\int_{b_1}^{b_2} g(\tau )e^{itf(\tau )}d\tau =
\sqrt{\frac{2\pi}{t|f^{\prime\prime}(\tau _{1})|}}g(\tau _{1})e^{itf(\tau _{1})+\frac{i\pi}{4}sgn f^{\prime\prime}(\tau _{1})} + O\left(\frac{g(b_1)}{t f'(b_1)}\right) \\ + O\left(\frac{g(b_2)}{t f'(b_2)}\right), \quad t \to \infty.
\end{multline}
Using \eqref{eq417} we obtain
\begin{equation*}
\frac{\partial ^{2}F(\tau _{1},\lambda)}{\partial \tau ^{2}}=\frac{1}{\tau _{1}(1-\tau _{1})}=\frac{(1+\lambda )^{2}}{\lambda }.
\end{equation*}
Evaluating $F(\tau,\lambda)$ at $\tau =\tau _{1}$ we find
\begin{multline*}
F(\tau _{1},\lambda)=\frac{\lambda }{1+\lambda }\ln \left(\frac{\lambda }{1+\lambda }\right)
+\frac{1 }{1+\lambda }\ln \left(\frac{1 }{1+\lambda }\right)
+\frac{1 }{1+\lambda }\ln\lambda \\
=\frac{\lambda }{1+\lambda }\ln \left(\frac{\lambda }{1+\lambda }\right)
+\frac{1 }{1+\lambda }\ln \left(\frac{1 }{1+\lambda }\right) = \ln{\frac{\lambda }{1+\lambda }} =  -\ln \left(1+\frac{1}{\lambda }\right).
\end{multline*}
Thus,
\begin{equation}  \label{eq420}
\sqrt{\frac{2\pi}{t|f^{\prime\prime}(\tau _{1})|}}e^{itf(\tau _{1})+\frac{i\pi}{4}sgn f^{\prime\prime}(\tau _{1})}
= \sqrt{\frac{2\pi}{t}}e^{\frac{i\pi}{4}} \frac{\lambda^{\frac{1}{2}}}{(1+\lambda)}(1+\lambda)^{-it} =\sqrt{\frac{2\pi}{t}}e^{\frac{i\pi}{4}}\frac{\lambda ^{\frac{1}{2}+it}}{(1+\lambda )^{1+it}}.
\end{equation}
The definition of $G(\sigma ,\tau )$ in the first of equations \eqref{eq44} implies
\begin{equation} \label{eq421}
g(\tau_1) = G(\sigma ,\tau _{1})=\bigl( 1-\frac{1}{1+\lambda } \bigr)^{-\frac{1}{2}}\frac{1}{(1+\lambda )^{\sigma -\frac{1}{2}}}=
\frac{\lambda ^{-\frac{1}{2}}}{(1+\lambda )^{\sigma -1}}.
\end{equation}
Substituting equations \eqref{eq411}, \eqref{eq420} and \eqref{eq421} into equation \eqref{eq419}, we find \eqref{eq415}.

Assuming that $1/(1+\lambda)$ does not approach $1-t^{\delta_3 - 1}$, it is straightforward to compute the large $t$-asymptotics of $\tilde{J}_3^U$ via integration by parts:
\begin{align} \label{eq422}
\tilde{J}_3^U &= \frac{1}{it} \int_{1-t^{\delta_3 -1}}^{\infty e^{i\varphi}} \frac{G}{\partial F / \partial \tau} \left( \frac{\partial }{\partial \tau} e^{itF} \right) d\tau \nonumber \\ &= - \frac{G e^{itF}}{it \partial F / \partial \tau} \Bigg|_{\tau=1-t^{\delta_3 -1}} - \frac{1}{it} \int_{1-t^{\delta_3 -1}}^{\infty e^{i\varphi}} \left[ \frac{\partial}{\partial \tau} \left( \frac{G}{\partial F / \partial \tau} \right) \right] e^{itF} d\tau.
\end{align}
Using the identities
\begin{equation*}
F(1-t^{\delta_3 -1}, \lambda) = t^{\delta_3 -1} \ln{t^{\delta_3 -1}} + (1-t^{\delta_3 -1}) \ln{(1-t^{\delta_3 -1})} + (1-t^{\delta_3 -1}) \ln{\lambda},
\end{equation*}
\begin{equation*}
\frac{\partial F}{\partial \tau}(1-t^{\delta_3 -1}, \lambda) = \ln{\left( \frac{\lambda (1-t^{\delta_3 -1})}{t^{\delta_3 -1}}  \right)},
\end{equation*}
together with equation \eqref{eq410} we find that the first term of the rhs of \eqref{eq422} yields the leading term of the rhs of \eqref{eq416}. The rigorous derivation of the relevant error term, as well as the analysis of the case that the stationary point approaches $1-t^{\delta_3 -1}$, is presented in [FSF].
\textbf{QED}

\begin{remark} \label{rem41}
The contour $L(\delta_2)$ can be deformed to a contour which can be written in  the form 
\begin{equation*}
L(\delta_2)=\left[ t^{\delta_2-1}, -\infty e^{i\Phi}\right]\cup \left[ -\infty e^{i\Phi},\infty e^{i\varphi}\right],
\end{equation*}
where $\Phi$ is appropriately constrained so that the associated integral converges. Hence, $J_3^S$ can be written
\begin{equation*}
J_3^S=J_3^{SP}+J_3^L
\end{equation*}
The leading behaviour of  $J_3^{SP}$ is given by the rhs of \eqref{eq415} and the leading behaviour of  $J_3^{L}$ is given via integration by parts, assuming that the endpoint $ t^{\delta_2-1}$ does not approach a stationary point:
\begin{equation}\label{422a}
\tilde{J}_3^L (\sigma,t,\delta_2, \lambda) = - \frac{1}{it} \frac{G e^{itF}}{\frac{\partial F}{\partial \tau}} \Bigg|_{\tau=t^{\delta_2 -1}}.
\end{equation}
Using the identities
\begin{equation*}
F(t^{\delta_2 -1}, \lambda) = (1-t^{\delta_2 -1}) \ln{(1-t^{\delta_2 -1})} + t^{\delta_2 -1} \ln{t^{\delta_2 -1}} + t^{\delta_2 -1} \ln{\lambda},
\end{equation*}
\begin{equation*}
\frac{\partial F}{\partial \tau}(t^{\delta_2 -1}, \lambda) = \ln{\left( \frac{\lambda}{t^{1-\delta_2 }-1}  \right)},
\end{equation*}
together with equation \eqref{eq411} we find
\begin{equation} \label{eq423}
\tilde{J}_3^L (\sigma,t,\delta_2, \lambda) =  \frac{it^{\delta_2 \left( \sigma - \frac{1}{2}  \right)}}{t^{\sigma + \frac{1}{2} }} \frac{(1-t^{\delta_2 -1})^{-\frac{1}{2}+i(t-t^{\delta_2})}}{\ln{\left(  \frac{\lambda}{t^{1-\delta_2}-1} \right)}} \lambda^{it^{\delta_2}} t^{i(\delta_2 -1)t^{\delta_2}}.
\end{equation}
The analysis of Remark 5.2 indicates that the relevant contribution is negligible. The rigorous derivation of \eqref{422a} as well as the analysis of the case that the stationary point approaches the lower end point $t^{\delta_2 -1}$ is similar with the analysis presented in \cite{FSF}; details are given in \cite{F}.
\end{remark}

\begin{remark}
Employing \eqref{eq423} in $\tilde{I}_3$  with $\sigma=1/2$ yields the leading contribution
\begin{align}
- &\frac{1}{\sqrt{2\pi}} \dfrac{1}{\sqrt{t-t^{\delta_2}}} e^{\frac{i\pi}{4}}t^{i(\delta_2-1)t^{\delta_{2}}}   \left(1-t^{\delta_2-1}\right)^{i(t-t^{\delta_2})}\\
& \times \mathop{\sum\sum}_{m_1,m_2\in \tilde{N}(\delta_2,t)}\frac{1}{m_{1}^{\frac{1}{2}+it^{\delta_2}}m_{2}^{\frac{1}{2}-it^{\delta_2}}} \dfrac{1}{\ln\left[\frac{m_1}{m_2}\left(t^{1-\delta_2}-1\right)\right]}, \qquad 0<\delta_{2}<1. \notag
\end{align}
where $\tilde{N}(\delta_2,t)=M_t\cap \tilde{M}_r^c$, with
\begin{equation*}
 M_t=\Big\{m_{j}=1,\ldots,[t], \ j=1,2\Big\} 
 \end{equation*}
and  $\tilde{M}_r^c$ denotes the complement of $\tilde{M}_r$, which is given by
\begin{equation*}\tilde{M}_r(\delta_2,t)=\left\{(m_1,m_2), \ \frac{m_{1}}{m_{2}}=t^{\delta_{2}-1}\left(1+ O\left(t^{-\epsilon}\right)\right), \ \epsilon>0\right\},
 \end{equation*}
 with $\epsilon>0$ is arbitrarily small.

The following heuristic argument implies that for the interesting case of  $\delta_2=1/2$, the contribution of the error term is negligible: employing (1.3) of \cite{FL} with $\eta\to2\pi t, \ t\to t^{\delta}, \ \sigma=1/2, \ \delta\in(0,1)$, we obtain
\begin{equation*}
\zeta\left(\frac{1}{2}+i t^\delta\right)=\sum_{n=1}^{[t]} \dfrac{1}{n^{\frac{1}{2}+i t^\delta}} + O\left( t^{\frac{1}{2}-\delta}\right), \quad t \to \infty,
\end{equation*}
thus
\begin{align*}
\sum_{m_{1}=1}^{[t]} \sum_{m_{2}=1}^{[t]}\frac{1}{m_{1}^{\frac{1}{2}+it^{\delta}}m_{2}^{\frac{1}{2}-it^{\delta}}} =\left|\zeta\left(\frac{1}{2}+i t^\delta\right)\right|^2 &+ O\left( t^{\frac{1}{2}-\delta}\left|\zeta\left(\frac{1}{2}+i t^\delta\right)\right|\right)\\
& +  O\left( t^{1-2\delta}\right),\quad t \to \infty.
\end{align*}
Using the estimate $\zeta\left(\frac{1}{2}+i t^\delta\right)=O\left( t^{\frac{\delta}{6}}\right),$ the above expression yields
\begin{align*}
\frac{1}{\sqrt{t}}\sum_{m_{1}=1}^{[t]} \sum_{m_{2}=1}^{[t]}\frac{1}{m_{1}^{\frac{1}{2}+it^{\delta}}m_{2}^{\frac{1}{2}-it^{\delta}}} =O\left( t^{\frac{\delta}{3}-\frac{1}{2}}\right) + O\left( t^{-\frac{5}{6}\delta}\right) +  O\left( t^{\frac{1}{2}-2\delta}\right),\quad t \to \infty.
\end{align*}
The last equation suggests that for $\delta_2>1/4$ the steepest descent contribution of \emph{(5.25)} is bounded by a decreasing
function of $t$, which vanishes as $t\to\infty$.

Rigorous estimates of the steepest descent contribution are presented in \cite{F}.
\end{remark}

\begin{theorem} \label{theo41}
Let $\tilde{I}_3$ be defined in \eqref{1.10}. Then,
\begin{multline} \label{eq424}
\tilde{I}_3 (\sigma,t, \delta_2, \delta_3)  = \\ 2 \Re \left\{ \mathop{\sum\sum}_{m_1,m_2 \in M(\delta_2, \delta_3)} \frac{1}{m_2^{\bar{s}}(m_1+m_2)^s} %+ \frac{1}{m_1^{\sigma} m_2^{\sigma}} O \left(  \frac{t^{\delta_2 \left( \sigma - \frac{1}{2} \right)}}{t^{\sigma}}, \lambda \right)
[1+o(1)]\right\} \left[ 1 + O(t^{-\delta_{23}})  \right] \\ 
- \sqrt{\frac{2}{\pi}} \Re {\left\{ e^{\frac{i\pi}{4}} t^{-\frac{\delta_3}{2}} (1-t^{\delta_3 -1})^{\sigma - \frac{1}{2} +i(t-t^{\delta_3})} t^{i(\delta_3-1)t^{\delta_3}} %\right. }\\ \left.
\mathop{\sum\sum}_{m_1,m_2 \in N(\delta_3,t)} \frac{1}{m_1^{s-it^{\delta_{3}}}} \frac{1}{m_2^{\bar{s}+it^{\delta_{3}}}} \right. }\\ \left. \times \frac{1}{\ln{\left[ \frac{m_2}{m_1} \left( t^{1-\delta_{3}} - 1 \right)  \right]}} [1+o(1)]  %+ O(t^{-\delta_3}, \lambda)
\right\} \left[  1 + O(t^{-\delta_{23}}) \right]  \\
+ \Re \left\{ \mathop{\sum\sum}_{m_1,m_2 \in M_r(\delta_3,t)} \frac{1}{m_1^s} \frac{1}{m_2^{\bar{s}}} E_3^T\left(\delta_3,t,\frac{m_2}{m_1})\right)\right\}\left[  1 + O(t^{-\delta_{23}}) \right], \quad t \to\infty,
\end{multline}
where $\sigma,\delta_2,\delta_3,\delta_{23}$, as well as the set $M$ is given in \eqref{1.14}, the sets $N$ and $M_r$ are as in and \emph{(1.24)}, and $E_3^T$ denotes the contribution of the points belonging in the transition zone.
\end{theorem}

\noindent \textbf{Proof} Expressing in \eqref{1.10} $\tilde{J}_3$ in terms of $\tilde{J}_3=\tilde{J}_3^S - \tilde{J}_3^U$, using equations \eqref{eq415} and  \eqref{eq416}  for $\tilde{J}_3^S$ and $\tilde{J}_3^U$ respectively, and simplifying the resulting formulae, we find \eqref{eq424}. \textbf{QED}

\section{The Leading asymptotics of $I_4$}

Let $I_4$ be defined by
\begin{equation} \label{5.1}
I_4(\sigma, t,\delta_2,\delta_3) = \frac{t}{\pi} \oint_{1-t^{\delta_3-1}}^{1+t^{\delta_4-1}} \Re\left\{ \frac{\Gamma(it - it\tau)  }{\Gamma(\sigma + i t)}\Gamma(\sigma + it\tau)\right\} | \zeta(\sigma+it\tau)|^2 \textrm{d}\tau,   \nonumber
\end{equation}
\begin{equation}
\quad 0<\sigma<1, \quad t>0,\quad 0<\delta_{2}<1, \quad 0<\delta_{3}<1,
\end{equation}
where the principal value integral is with respect to $\tau =1$. Using the change of variables $1-\tau=\rho$, $I_4$ becomes
\begin{equation*}
I_4 = \frac{t}{\pi} \oint_{-t^{\delta_4-1}}^{t^{\delta_3-1}} \Re{\left\{ \Gamma(it\rho)  \frac{\Gamma(\sigma + i t -it\rho)}{\Gamma(\sigma + i t)}\right\}} \left| \zeta(\sigma+i t-it\rho)\right|^2 \textrm{d}\rho,
\end{equation*}
where now the principal value integral is defined with respect to $\rho=0$. The change of variables $t\rho=x$ yields
\begin{equation} \label{5.4}
I_4 = \frac{1}{\pi} \oint_{-t^{\delta_4}}^{t^{\delta_3}} \Re{\left\{ \Gamma(ix)  \frac{\Gamma(\sigma + i t -ix)}{\Gamma(\sigma + i t)}\right\}} \left| \zeta(\sigma+i t-ix)\right|^2 \textrm{d}x.
\end{equation}
In the interval of integration we have
\begin{equation*}
\quad -t^{\delta_4}\leq x\leq t^{\delta_3}.
\end{equation*}
Thus,
\begin{equation}
\quad t-t^{\delta_3}\leq t-x \leq t+t^{\delta_4}.
\end{equation}

Replacing in (6.2), $|\zeta|^{2}$ by its leading order asymptotics we find
\begin{equation}
\zeta\left(\sigma + i(t-x)\right) \sim  \sum_{m=1}^{[\frac{\eta}{2\pi}]} \frac{1}{m^{\sigma+i(t-x)}}, \quad \eta > t-x.
\end{equation}
Since $t-x \leq t + t^{\delta 4}$, we take $\eta = 2\pi t > t + t ^{\delta_4}$. Thus,
\begin{equation} \label{5.13}
|\zeta(\sigma +i(t - x)) |^2 \sim  \sum_{m_1=1}^{[t]}  \sum_{m_2=1}^{[t]}  \frac{1}{m_1^s m_2^{\bar{s}}} \left(\frac{m_1}{m_2}\right)^{ix}.
\end{equation}

Let $\tilde{I_4}$ denote the expression obtained from $I_{4}$ by replacing $|\zeta|^{2}$ with the rhs of (6.5), i.e,
%Substituting \eqref{5.9} and \eqref{5.13} into \eqref{5.6} we find
%\begin{align} \label{5.14}
%I_4(\sigma, t,\delta_3,\delta_4) &= I_4^{(1)}(\sigma, t,\delta_3,\delta_4) + \overline{G(e^{it})}T^{it} \frac{T^{-\sigma}}{\pi} I_4^{(2)}(\sigma, t,\delta_3,\delta_4) \nonumber \\
%& \quad+ G(e^{it}) T^{-it} \frac{T^{-\sigma}}{\pi} I_4^{(3)}(\sigma, t,\delta_3,\delta_4) + |G|^2 \frac{T^{-2\sigma}}{\pi^2} I_4^{(4)}(\sigma, t,\delta_3,\delta_4) \nonumber \\
%& \quad+  I_4^{(5)}(\sigma, t,\delta_3,\delta_4) + I_4^{(6)}(\sigma, t,\delta_3,\delta_4),
%\end{align}
%where
\begin{align}  \label{5.15a}
&\tilde{I_4}(\sigma, t,\delta_3,\delta_4) =\Re\left\{  \sum_{m_1=1}^{[t]}  \sum_{m_2=1}^{[t]} \frac{1}{m_{1}^s}\frac{1}{m_{2}^{\bar{s}}} J_4(\sigma, t,\delta_{2},\delta_{3},\frac{m_{1}}{m_{2}})\right\},
\end{align}
%\begin{multline} \label{5.15b}
%I_4^{(2)}(\sigma, t,\delta_3,\delta_4) = \frac{1}{\pi} \Re{\left\{\sum_{m=1}^{\left \lfloor{T}\right \rfloor} \frac{1}{m^s} \int_{H} \frac{e^z}{z}J_4(\sigma,t,\delta_3,\delta_4,A_2)\textrm{d}z \right\}} \\
%\times \left[ 1 + O\left( \frac{1}{t} \right) \right], \quad t\to\infty,
%\end{multline}
%\begin{multline} \label{5.15c}
%I_4^{(3)}(\sigma, t,\delta_3,\delta_4) = \frac{1}{\pi} \Re{\left\{\sum_{m=1}^{\left \lfloor{T}\right \rfloor} \frac{1}{m^{\bar{s}}} \int_{H} \frac{e^z}{z}J_4(\sigma,t,\delta_3,\delta_4,A_3)\textrm{d}z \right\}} \\
%\times \left[ 1 + O\left( \frac{1}{t} \right) \right], \quad t\to\infty,
%\end{multline}
%\begin{multline} \label{5.15d}
%I_4^{(4)}(\sigma, t,\delta_3,\delta_4) = \frac{1}{\pi} \Re{\left\{ \int_{H} \frac{e^z}{z}J_4(\sigma,t,\delta_3,\delta_4,A_4)\textrm{d}z \right\}} \\
%\times \left[ 1 + O\left( \frac{1}{t} \right) \right], \quad t\to\infty,
%\end{multline}
%\begin{equation} \label{5.15e}
%I_4^{(5)}(\sigma, t,\delta_3,\delta_4) = \frac{1}{\pi} \Re{\left\{\sum_{m=1}^{\left \lfloor{T}\right \rfloor} \frac{1}{m^s} \int_{H} \frac{e^z}{z}\hat{J}_4(\sigma,t,\delta_3,\delta_4,A_5)\textrm{d}z \right\}},
%\end{equation}
%\begin{equation} \label{5.15f}
%I_4^{(6)}(\sigma, t,\delta_3,\delta_4) = \frac{1}{\pi} \Re{\left\{\sum_{m=1}^{\left \lfloor{T}\right \rfloor} \frac{1}{m^{\bar{s}}} \int_{H} \frac{e^z}{z}\hat{J}_4(\sigma,t,\delta_3,\delta_4,A_6)\textrm{d}z \right\}},
%\end{equation}
%\end{subequations}
where $J_4$ is defined by
\begin{align} \label{5.17}
&J_4(\sigma,t,\delta_3,\delta_4,\frac{m_{1}}{m_{2}}) =\frac{1}{\pi}\oint_{-t^{\delta_4}}^{t^{\delta_3}} \Gamma (i x) \frac{\Gamma (\sigma + i t - i x)}{\Gamma(\sigma + i t )}\left(\frac{m_{1}}{m_{2}}\right)^{i x }\textrm{d}x, \nonumber\\
&0<\sigma<1,~~t>0,~~0<\delta_{3}<1,~~0<\delta_{4}<1,~~m_{j}=1,2,\ldots, [t],
\end{align}
with the principal value integral  defined with respect to $x=0$.
\begin{proposition}
Let $J_{4}$ be defined by (6.7).  Let $H_{1}$ denote the Hankel contour with a branch cut along the negative real axis, see figure 1, defined by
\begin{align}
H_{1}=\left\{r e^{- i \pi}|1<r<\infty\right\}\cup\left\{e^{ i \theta}|-\pi<\theta<\pi\right\}\cup\left\{r e^{ i \pi}|1<r<\infty\right\}.
\end{align}
Then,
\begin{align}
&J_4(\sigma, t,\delta_3,\delta_4,\frac{m_{1}}{m_{2}})=\frac{1}{\pi}\int_{H_{1}}\frac{e^{z}}
{z}
\tilde{J_4}(\sigma, t,\delta_3,\delta_4,A)\textrm{d}z \left[1+O\left(\frac{1}{t}\right)\right],\ \ \ t\to\infty,\notag\\
&0<\sigma<1,~0<\delta_{3}<\frac{1}{2},~0<\delta_{4}<\frac{1}{2},~
A=\frac{m_1}{m_2}\frac{z}{t},
\end{align}
where
\begin{align}
\tilde{J_4}(\sigma, t,\delta_3,\delta_4, A) = \oint_{-t^{\delta_4}}^{t^{\delta_3}} \frac{e^{\frac{\pi x}{2}}A^{ix}}{e^{-\pi x}-e^{\pi x}} \left(1-\frac{x}{t}\right)^{\sigma-\frac{1}{2}} e^{i x}\left(1-\frac{x}{t}\right)^{i (t-x)}\textrm{d}x,
\end{align}
with the principal value integral  defined with respect to $x=0$.
\end{proposition}

\textbf{Proof}.
Equation (2.33a) together with the inequality
\begin{equation}
\frac{1}{t-x}\leq \frac{1}{t-t^{\delta_{3}}}=\frac{1}{t}+O\left(t^{\delta_{3}-2}\right),~~t\to\infty,
\nonumber
\end{equation}
imply
\begin{align}
\Gamma (\sigma + i (t-x))=& \sqrt{2\pi} (t-x)^{\sigma-\frac{1}{2}}e^{-\frac{\pi}{2}(t-x)}e^{-\frac{i\pi}{4}}e^{\frac{i\pi\sigma}{2}}\notag\\ 
&\times e^{-i(t-x)} (t-x)^{i(t-x)}\left(1+O\left(\frac{1}{t}\right)\right), \quad t\to\infty.
\end{align}
The above equation together with \eqref{eq34a} yield
\begin{equation}
\hspace*{-5mm}\frac{\Gamma \left(\sigma + i (t-x)\right)}{\Gamma (\sigma + it)}= \left(1-\frac{x}{t}\right)^{\sigma -\frac{1}{2}}e^{\frac{\pi x}{2}}t^{-ix}e^{ix}\left(1-\frac{x}{t}\right)^{i(t-x)}\left[1+O\left(\frac{1}{t}\right)\right],~t\to\infty.
\end{equation}
 Replacing in equation (6.7), $\Gamma (\sigma + i (t-x))/\Gamma (\sigma + i t)$ by the rhs of (6.12), as well as employing the formula
 \begin{equation}
\Gamma ( i x)= \frac{1}{e^{-\pi x}-e^{\pi x}} \int_{H_1}\frac{e^{z}}{z}z^{ix}\textrm{d}z,
\end{equation}
equation (6.7) becomes equation (6.9). \textbf{QED}

\begin{proposition}
Let $\tilde{J_4}$ be defined in (6.10).  Then,
\begin{align}
&  \tilde{J_4}(\sigma, t,\delta_3,\delta_4,A) =\left[\frac{i}{2}(-1+\frac{2}{1-i A})+\frac{e^{it^{\delta_{3}}\ln A} e^{- \frac{\pi t^{\delta_3}}{2}}}{\frac{\pi}{2}-i \ln A}\right] \left[1 +O\left(t^{-\delta_{34}}\right)\right], \nonumber\\
&t \to\infty, \ 0<\sigma <1,~0<\delta_{3}<\frac{1}{2},~0<\delta_{4}<\frac{1}{2},~\delta_{34}=\min\{\delta_{3},\delta_{4}\},
~A=\frac{m_{1}z}{m_{1}t}.
\end{align}
\end{proposition}
\textbf{Proof} Using the identity
\begin{align*} \label{5.19}
e^{ix}\left(  1 - \frac{x}{t} \right)^{i(t-x)} &= e^{ix} e^{i(t-x)\ln{\left(  1 - \frac{x}{t} \right)}} =e^{ix} e^{i(t-x) \left(  -\frac{x}{t} + O \left(  \frac{x^2}{t^2} \right) \right) } \\
&= e^{iO \left(  \frac{x^2}{t} \right) } = 1 + O \left(  \frac{x^2}{t} \right), \quad \frac{x}{t} \to 0,
\end{align*}
we find
\begin{equation}\nonumber
\left(1-\frac{x}{t} \right)^{\sigma - \frac{1}{2}}e^{ix} \left(  1 - \frac{x}{t} \right)^{i(t-x)} =\left[ 1 + O\left( t^{2\delta_{34}-1} \right) \right],~t\to\infty,
\end{equation}
where we have used that $|x|< t ^{\delta_{34}}$. Then, equation (6.10) becomes
\begin{equation}
\tilde{J}_4(\sigma, t,\delta_3,\delta_4,A) = \oint_{-t^{\delta_{4}}}^{t^{\delta_{3}}}\frac{e^{\frac{\pi x}{2}}A^{ix}}{e^{-\pi x}-e^{\pi x}} \textrm{d}x \left[1+O\left(t^{2\delta_{34}-1}\right)\right],~~t \rightarrow \infty.
\end{equation}

It is remarkable  that the leading order term of the above integral  can be computed in closed form within an error which is exponentially small as $t\to\infty$. Indeed,
\begin{equation}
 \lim_{\varepsilon\to 0} \left\{\int_{-t^{\delta_4}}^{-\varepsilon} \frac{e^{\frac{\pi x}{2}}A^{ix}}{e^{-\pi x}-e^{\pi x}} \textrm{d}x + \int_{\varepsilon}^{t^{\delta_3}} \frac{e^{\frac{\pi x}{2}}A^{ix}}{e^{-\pi x}-e^{\pi x}} \textrm{d}x \right\}
\nonumber
\end{equation}
\begin{equation}
= \lim_{\varepsilon\to 0} \left\{ \int_{\varepsilon}^{t^{\delta_4}} \frac{e^{-\frac{3\pi x}{2}}A^{-ix}}{1-e^{-2\pi x}} \textrm{d}x - \int_{\varepsilon}^{t^{\delta_3}} \frac{e^{-\frac{\pi x}{2}}A^{ix}}{1-e^{-2\pi x}} \textrm{d}x \right\}
\nonumber
\end{equation}
\begin{equation}
= \lim_{\varepsilon\to 0} \Biggl( \int_{\varepsilon}^{t^{\delta_4}} e^{-\frac{3\pi x}{2}}A^{-ix} \left( \sum_{k=0}^{\infty}{e^{-2\pi kx}} \right) \textrm{d}x - \int_{\varepsilon}^{t^{\delta_3}} e^{-\frac{\pi x}{2}}A^{ix} \left( \sum_{k=0}^{\infty}{e^{-2\pi kx}} \right) \textrm{d}x
\Biggr)
\nonumber
\end{equation}
\begin{equation}
= \lim_{\varepsilon\to 0} \sum_{k=0}^{\infty} \left\{ \int_{\varepsilon}^{t^{\delta_4}} e^{-x\left(2\pi k + \frac{3\pi}{2} + i\ln{A}\right)} \textrm{d}x - \int_{\varepsilon}^{t^{\delta_3}}  e^{-x\left(2\pi k + \frac{\pi}{2} - i\ln{A}\right)} \textrm{d}x\right\}
\nonumber
\end{equation}
\begin{equation}
=- \sum_{k=0}^{\infty} \left\{ \frac{e^{- t^{\delta_4} \left(2\pi k + \frac{3\pi}{2} + i\ln{A}\right)}-1}{2\pi k + \frac{3\pi}{2}+i\ln{A}} -  \frac{e^{- t^{\delta_3} \left(2\pi k + \frac{\pi}{2} - i\ln{A}\right)}-1}{2\pi k + \frac{\pi}{2}-i\ln{A}} \right\}.\nonumber
\end{equation}
The validity of the above interchange of the limit $\epsilon\to 0$ and  the sum over $k$  has to be treated carefully: standard methods, such as dominated convergence fail; it is shown in \cite{F} that the validity of this  interchange can be proven rigorously via the Vitali convergence theorem.

Taking into consistent that $\arg z \in[-\pi,\pi]$, it follows that the terms involving $t^{\delta_{3}}$ and  $t^{\delta_{4}}$ decay exponentially expect for the term involving  $t^{\delta_{3}}$ and $k=0$. Hence,
\begin{align} \label{5.22}
\oint_{-t^{\delta_{4}}}^{t^{\delta_{3}}}\frac{e^{\frac{\pi x}{2}}A^{ix}}{e^{-\pi x}-e^{\pi x}} \textrm{d}x =& \sum_{k=0}^{\infty} \left\{ \frac{1}{2\pi k + \frac{3\pi}{2}+i\ln{A}} - \frac{1}{2\pi k + \frac{\pi}{2}-i\ln{A}}  \right\}\notag \\
&+\frac{e^{it^{\delta_{3}}\ln A} e^{- \frac{\pi t^{\delta_{3}}}{2}}}{\frac{\pi}{2}-i \ln A}+ O\left( e^{-t^{\delta_{34}}} \right), \quad t\to\infty.
\end{align}
Let $S$ denote the first term of the rhs of (6.16). Then,
%\begin{equation} \label{5.24}
%S= -\frac{\pi+2i\ln{A}}{(2\pi)^2} \sum_{k=0}^{\infty}  \frac{1}{(k+a)(k+b)} ,
%\end{equation}
\begin{equation} \label{5.24}
S= \frac{1}{2\pi}\sum_{k=0}^\infty\left(\frac{1}{k+1-b}-\frac{1}{k+b}\right) ,
\end{equation}
with
%where $a$ and $b$ are defined by
\begin{equation} \label{5.25}
%a=1-b, \qquad 
b=\frac{1}{4} - \frac{i}{2\pi}\ln{A}.
\end{equation}
%In order to simplify $\tilde{J}_4$ we will employ the identity
%\begin{equation} \label{5.19}
%\sum_{k=0}^{\infty}  \frac{1}{(k+a)(k+b)} = \frac{\Psi(a)-\Psi(b)}{a-b},\tag{5.19}
%\end{equation}
%where $\Psi(z)$, $z\in\mathbb{C}$, denotes the digamma function defined in (2.11). Employing in (6.19) the reflection formula for  $\Psi(z)$, namely
%\begin{equation} \label{5.20}
%\Psi(1-z) -\Psi(z) = \pi\cot{(\pi z)},\tag{5.20}
%\end{equation}
%we find
%\begin{equation} \label{5.21}
%\sum_{k=0}^{\infty}  \frac{1}{(k+1-b)(k+b)} = \frac{ \pi\cot{(\pi b)}}{1-2b}.\tag{5.21}
%\end{equation}
%Using in the rhs of (6.21) the definition of $b$ given in the second of equations (6.18) we find the following:
%\begin{equation*}
%\frac{ \pi\cot{(\pi b)}}{1-2b} = \frac{ \pi\cot{\left(\frac{\pi}{4} - \frac{i}{2} \ln{A}\right)}}{1-2\left(\frac{1}{4} - \frac{i}{2\pi} \ln{A}\right)}  = - \frac{ \pi\cot{\left(\frac{i}{2} \ln{A} +\frac{\pi}{4} - \frac{\pi}{2} \right)}}{\frac{1}{2} + \frac{i}{\pi} \ln{A}}
%\nonumber
%\end{equation*}
%\begin{equation*}
%=  \frac{ \pi\tan{\left(\frac{\pi}{4} + \frac{i}{2} \ln{A} \right)}}{\frac{1}{2\pi} (\pi+ 2i\ln{A})}.
%\end{equation*}
%Substituting this expression in (6.21) we find
%\begin{equation*} \label{5.29}
%S = -\frac{1}{2}\tan\left(\frac{\pi}{4} + \frac{i}{2} \ln{A} \right).
%\end{equation*}

Let $\Psi(z)$, $z\in\mathbb{C}$, denotes the digamma function defined in (2.11), then (6.17) takes the form 
\begin{equation} \label{5.19}
S=-\frac{1}{2\pi}\left(\Psi(1-b)-\Psi(b)\right).
\end{equation}
Employing in (6.19) the reflection formula for  $\Psi(z)$, namely
\begin{equation} \label{5.20}
\Psi(1-z) -\Psi(z) = \pi\cot{(\pi z)},
\end{equation}
we find
\begin{equation} \label{5.29}
S = -\frac{1}{2}\tan\left(\frac{\pi}{4} + \frac{i}{2} \ln{A} \right).
\end{equation}
This formula can be further simplified as follows:
\begin{equation*}
\tan{\left(\frac{\pi}{4} + \frac{i}{2} \ln{A} \right)} = \frac{1}{i} \frac{e^{i\left(\frac{\pi}{4} + \frac{i}{2} \ln{A}\right)} - e^{-i\left(\frac{\pi}{4} + \frac{i}{2} \ln{A}\right)}} {e^{i\left(\frac{\pi}{4} + \frac{i}{2} \ln{A}\right)} + e^{-i\left(\frac{\pi}{4} + \frac{i}{2} \ln{A}\right)}}
\nonumber
\end{equation*}
\begin{equation*}
= \frac{1}{i} \frac{e^{\frac{i\pi}{4}}A^{-\frac{1}{2}} - e^{\frac{-i\pi}{4}}A^{\frac{1}{2}}} {e^{\frac{i\pi}{4}}A^{-\frac{1}{2}} + e^{\frac{-i\pi}{4}}A^{\frac{1}{2}}} = \frac{1+iA} {i(1-iA)} = \frac{1}{i} \left( -1 + \frac{2}{1-iA} \right).
\end{equation*}
Hence
\begin{equation} \label{5.22a}
S = \frac{i}{2} \left( -1 + \frac{2}{1-iA} \right).
\end{equation}

Replacing in (6.15) the leading order term with the rhs of (6.16), where $S$ is given by (6.22), we find (6.14).
\textbf{QED}

\begin{proposition} Let $J_{4}$ be defined by (6.7). Then,
\begin{equation*} \label{5.23a}
J_{4}(\sigma, t, \delta_{3}, \delta_{4})=\left[-1+E_{4}(t, \delta_{3},M)\right][1+O(t^{2\delta_{34}-1})],\quad t\to\infty, \tag{6.23}
\end{equation*}
where
\begin{equation*} \label{5.24a}
E_{4}( t, \delta_{3}, M)=\frac{1}{\pi}\int_{H_{1}}\frac{e^{z}}{z} \frac{e^{it^{\delta_{3}}\ln(Mz)-\frac{\pi  t^{\delta_{3}}}{2}}}{\frac{\pi}{2}-i\ln(Mz)}dz,\ \ \ M=\frac{m_{1}}{m_{2}t}.\tag{6.24}
\end{equation*}
\end{proposition}
\textbf{Proof} Equation (6.23) follows from equations (6.9) and (6.14) with the aid of the following identity:
\begin{equation*} \label{5.25a}
\frac{i}{2\pi}\int_{H_{1}}\frac{e^{z}}{z}\left(-1+\frac{2}{1-iA}\right)dz
=-1,\ \ \ A=\frac{m_{1}}{m_{2}}\frac{z}{t}.\tag{6.25}
\end{equation*}
In order to derive (6.25) we will employ the following residue formulae:
\begin{equation} \nonumber
\int_{H_{1}} \frac{e^z}{z}\textrm{d}\tau = 2\pi i,
\end{equation}
and
\begin{equation*} \label{5.32a}
\int_{H_1} \frac{e^z}{z(z+ic)}\textrm{d}z = 2\pi i  \underset{z=0}{\operatorname{Res}}{\frac{e^z}{z(z+ic)}} = \frac{2\pi}{c}, \quad c\ne 0.
\end{equation*}
These formulae imply the following identities for the two terms occurring in the lhs of (6.25):
\begin{equation} \nonumber
-\frac{i}{2}\int_{H_1} \frac{e^z}{z}\textrm{d}z = \pi,
\end{equation}
and
\begin{equation*}
i\int_{H_1} \frac{e^z}{z(1-iA)}\textrm{d}z = \frac{i}{-i\left( \frac{m_1}{m_2t} \right)} \int_{H_1} \frac{e^z}{z\left(z+ \frac{im_2t}{m_1} \right)}\textrm{d}z = -2\pi.
\end{equation*}
Hence, (6.24) follows.

The pole of the integrand of $E_4$ occurs on the contour $H_1$ iff $m_1/m_2=t$. Letting $$m_2=1+\epsilon_1, \ \ m_1=t-\epsilon_2, \qquad \epsilon_1>0, \ \epsilon_2>0,$$
it follows that $$M\sim 1-\epsilon, \qquad \epsilon=\epsilon_1+ \frac{\epsilon_2}{t},  \quad \epsilon\to 0.$$
By employing the Plemelj formula it is straightforward to compute the limit $E_4$ as $\epsilon\to 0$. Details will be presented in \cite{F}.
\textbf{QED}

\begin{proposition} Let $E_{4}$ be defined in (6.24). Then,
\begin{equation*} \label{5.26}
E_{4}(t, \delta_{3}, M)=2e^{-\frac{i}{M}}+E_{4}^{SD}(t,\delta_{3}, M), \ \ \ M=\frac{m_{1}}{m_{2}t},\tag{6.26}
\end{equation*}
where the first term occurs iff
\begin{equation*} \label{5.27}
\frac{m_{1}}{m_{2}}\in(t^{1-\delta_{3}}, t),\tag{6.27}
\end{equation*}
and $E_{4}^{SD}$ is defined by
\begin{equation*}\label{5.28}
E_{4}^{SD}=\frac{1}{\pi}\int_{H_{1}}\frac{e^{t^{\delta_{3}}[\omega-\frac{\pi}{2}+i\ln(Mt^{\delta_{3}}\omega)] }      }{\omega[\frac{\pi}{2}-i\ln(Mt^{\delta_{3}}\omega)]}d\omega.\tag{6.28}
\end{equation*}
\end{proposition}
\textbf{Proof} In order to estimate $E_{4}$ we let $z=t^{\delta}\omega$. Then,
\begin{equation*}\label{5.29a}
E_{4}=\frac{1}{\pi}\int_{H_{t^{-\delta_{3}}}}\frac{e^{t^{\delta_{3}}[\omega-\frac{\pi}{2}+i\ln(Mt^{\delta_{3}}\omega)] }      }{\omega[\frac{\pi}{2}-i\ln(Mt^{\delta_{3}}\omega)]}d\omega,\tag{6.29}
\end{equation*}
where $H_{t^{-\delta_3}}$ is the Hankel contour involving a circle of radius  $t^{-\delta_{3}}$. The above integral has a stationary point at
\begin{equation*}\label{5.30}
\omega_{sp}=-i.\tag{6.30}
\end{equation*}
Thus, in order to estimate $E_4$, we deform the above circle in the $\omega$-complex plane to a circle of radius 1:
\begin{equation}\nonumber
E_{4}=\frac{1}{\pi}\int_{H_{1}}\frac{e^{t^{\delta_{3}}[\omega-\frac{\pi}{2}+i\ln(Mt^{\delta_3}\omega)]}      }{\omega[\frac{\pi}{2}-i\ln(Mt^{\delta_{3}}\omega)]}d\omega +E_{4}^{P},
\end{equation}
where $E_{4}^{P}$ is the contribution of the pole
\begin{equation*}\label{5.31}
\omega_{p}=-\frac{i}{t^{\delta_{3}}M}.\tag{6.31}
\end{equation*}
This pole contribution occurs, iff
\begin{equation*}\nonumber
\frac{1}{t^{\delta_{3}}M}\in(t^{-\delta_{3}},1),
\end{equation*}
or
\begin{equation*}
\frac{m_{2}t}{m_{1}t^{\delta_3}}\in(t^{-\delta_{3}},1),
\end{equation*}
which implies equation (6.27).
The residue of the integral of the rhs of (6.29) associate with $\omega_{p}$ is given by $ie^{-i/M}$. Cauchy's theorem yields  
\begin{equation}\label{5.32-new}
\frac{1}{\pi}\left(\int_{H_1}-\int_{H_{t^{-\delta_3}}}\right)\frac{e^{t^{\delta_{3}}[\omega-\frac{\pi}{2}+i\ln(Mt^{\delta_3}\omega)]}      }{\omega[\frac{\pi}{2}-i\ln(Mt^{\delta_{3}}\omega)]}d\omega = -2 e^{-i/M},\tag{6.32}
\end{equation}
and then equation (6.26) follows. 

The case that the pole $\omega_p$ approaches the stationary point $\omega_{sp}$ is analysed in \cite{F}.
\textbf{QED}

\begin{theorem} Let $\tilde{I}_{4}$ denote the integral obtained from $I_{4}$  defined in \eqref{1.7} with $j=4$, with $|\zeta|^{2}$  replaced by its leading term asymptotics. Then,
\begin{align}\label{5.32}
&\tilde{I}_{4}(\sigma, t, \delta_{3}, \delta_{4})=-\sum_{m_1=1}^{[t]}\sum_{m_{2}=1}^{[t]}
\frac{1}{m_{1}^{s}m_{2}^{\bar{s}}}\left[1+O(t^{2\delta_{34}-1})\right]\notag\\
&+2\Re\left\{\mathop{\sum\sum}_{m_1,m_2\in M_4(\delta_3,t)}
\frac{1}{m_{1}^{s}m_{2}^{\bar{s}}}e^{-\frac{im_{2}}{m_{1}}t}\right\}
\left[1+O(t^{2\delta_{34}-1})\right]\notag\\
&+\Re\left\{\sum_{m_1=1}^{[t]}\sum_{m_{2}=1}^{[t]}
\frac{1}{m_{1}^{s}m_{2}^{\bar{s}}}E_{4}^{SD}(t,\delta_{3},M)\right\}
\left[1+O(t^{2\delta_{34}-1})\right],~t\rightarrow\infty,\notag\\
&0<\sigma<1,~0<\delta_{3}<\frac{1}{2},~0<\delta_{4}<\frac{1}{2}, \ 
M=\frac{m_1}{m_2 t},
\tag{6.33}
\end{align}
where the set $M_4$ is defined by
\begin{equation*}
M_4(\delta_3,t)=\Big\{m_{j}=1,\ldots,[t],\ j=1,2,\
~\frac{m_{1}}{m_{2}}\in(t^{1-\delta_{3}},t)\Big\},
\end{equation*} 
and  $E_{4}^{SD}$ is defined by \emph{(6.28)} with $M=\frac{m_{1}}{m_{2}t}$. 
\end{theorem}
\textbf{Proof} $\tilde{I}_{4}$ can be expressed in terms of $J_{4}$ by equation (6.6) and $J_{4}$ is given by (6.23). Replacing in the latter equation $E_{4}$ by the rhs of (6.26) we find (6.33).
\textbf{QED}

\begin{remark}
 A steepest descent computation, for $(m_1,m_2)\in N(\delta_3)$, implies that the leading order term of $E_{4}^{SD}$ is given by
\begin{equation}
E_{4}^{SD}\sim-\sqrt{\frac{2}{\pi}}e^{\frac{i\pi}{4}} t^{-\frac{\delta_3}{2}}e^{-it^{\delta_3}} t^{i(\delta_3-1) t^{\delta_3}}\frac{1}{\ln\left(\frac{m_2}{m_1}t^{1-\delta_3}\right)} \left(\frac{m_1}{m_2}\right)^{it^{\delta_3}}, \ \ \ t\to\infty.  \tag{6.34}
\end{equation}

For completeness, $E_{4}^{SD}$ requires the analysis of the transition zone that now corresponds to the set of points where the pole approaches the steepest descent points, namely points  $(m_1,m_2)\in M_r(\delta_3)$; the rigorous computation is given in \cite{F}. 

%We note that the contribution from the whole pole yields $|E_4^P|=2$ and the steepest descent contribution yields $|E_4^{SD}|\sim \sqrt{\frac{2}{\pi}} \frac{1}{t^{\frac{\delta_3}{2}} \ln\left(\frac{m_2}{m_1}t^{1-\delta_3}\right)}$. Hence, by construction, $|E_4^T|\in\left(\sqrt{\frac{2}{\pi}} \frac{1}{t^{\frac{\delta_3}{2}} \ln\left(\frac{m_2}{m_1}t^{1-\delta_3}\right)},1\right).$
\end{remark}

%\begin{remark}\label{rem5}
%Recall that
%\begin{equation}\label{5.33}
%\left(1-\frac{x}{t}\right)^{\sigma-\frac{1}{2}}e^{ix}\left(1-\frac{x}{t}\right)^{i(t-x)}
%=\left[1+O\left(\frac{x}{t}\right)\right]\left[1+O\left(\frac{x^{2}}{t}\right)\right],
%~~t\rightarrow\infty.\tag{6.35}
%\end{equation}
%This estimate, together with the fact that $|x|<t^{\delta_{34}}$ implies the restriction $x^{2}/t<t^{2\delta_{34}-1}$, which then imposes that both $\delta_{3}$ and $\delta_{4}$ must be less than $1/2$. However, the constraint $x^{2}/t<t^{2\delta_{34}-1}$ can be eliminated and hence $\delta_{3}$ and $\delta_{4}$ only need to be constrained to be less than 1: one can estimate the integral $\tilde{J}_{4}$ defined in (6.10) directly without eliminating the term defined by the lhs of (6.35). 
%%The relevant approach is very similar to the one used in section 7, where it is no possible to eliminate the lhs of (6.35) since in that case $x^{2}/t$ may grow.
%
%
%\end{remark}

\section{Further Developments}

Define the integral $J$ by
\begin{align}\label{f.1}
J(\delta_1,\delta_2,t)=\sqrt{t}& \int_{t^{\delta_1-1}}^{t^{\delta_2-1}} g(t,\tau) e^{itf(\tau)} \left|\zeta\left(\frac{1}{2}+i\tau t\right)\right|^2 d\tau,\notag\\
& 0<\delta_j<1, \  j=1,2, \ \delta_2>\delta_1, \ t>0,
\end{align}
where $g$ and $f$ are given real functions. Replacing  $\left|\zeta\left(\frac{1}{2}+i t\right)\right|^2$ by its leading asymptotic sum $S_R$ defined in \eqref{asym_R} we find 
\begin{equation}\label{f.2}
J(\delta_1,\delta_2,t)\sim \sum_{m_1=1}^{[t]} \sum_{m_2=1}^{[t]} \frac{1}{\sqrt{m_1 m_2}} \tilde{J}(\delta_1,\delta_2,t,\lambda), \quad t\to\infty,
\end{equation}
where $\tilde{J}$ is defined by
\begin{equation}\label{f.3}
\tilde{J}(\delta_1,\delta_2,t,\lambda)=\sqrt{t}\int_{t^{\delta_1-1}}^{t^{\delta_2-1}} g(t,\tau) e^{itF(\tau,\lambda)}  d\tau,
\end{equation}
with 
\begin{equation}\label{f.4}
F(\tau,\lambda)=f(\tau)+\tau \ln \lambda, \qquad \lambda=\frac{m_2}{m_1}.
\end{equation}

\subsection*{Example 1}

Let 
\begin{equation}\label{f.5}
f(\tau)= (1-\tau) \ln(1-\tau) + \tau \ln \tau.
\end{equation}
Then,
\begin{align*}
F_\tau=-\ln(1-\tau) +\ln\tau + \ln\lambda,\\
F_{\tau\tau}=\frac{1}{\tau(\tau-\lambda)}.
\end{align*}
The stationary point denoted by $\xi$ satisfies
\begin{equation}\label{f.6}
\xi=\frac{1}{1+\lambda}, \qquad t^{\delta_1-1}\leq \frac{1}{1+\lambda} \leq  t^{\delta_2-1}.
\end{equation}
Using in equation (5.19) the relation $$F_{\tau\tau}\bigg|_{\tau=\xi}=\frac{(1+\lambda)^2}{\lambda},$$ and assuming that the main contribution to the large $t$-asymptotics of $\tilde{J}$ comes from $\tau=\xi$, we find
\begin{equation}\label{f.7}
\tilde{J}\sim \sqrt{2\pi} e^{i\frac{\pi}{4}} \frac{\sqrt{\lambda}}{1+\lambda} g\left(t,\frac{1}{1+\lambda}\right) e^{it F\left(\frac{1}{1+\lambda},\lambda\right)}. 
\end{equation}
But, \begin{align*}
F\left(\frac{1}{1+\lambda},\lambda\right)&=\left(1-\frac{1}{1+\lambda}\right) \ln\left(1-\frac{1}{1+\lambda}\right)+ \frac{1}{1+\lambda}\ln\left(\frac{1}{1+\lambda}\right) + \frac{1}{1+\lambda} \ln\lambda\\
&=\ln\frac{\lambda}{1+\lambda}.
\end{align*}
Thus, simplifying \eqref{f.7} and then substituting the resulting expression in \eqref{f.2} we find
\begin{equation}\label{f.8}
J(\delta_1,\delta_2,t)\sim \mathop{\sum\sum}_{m_1,m_2\in N(\delta_1,\delta_2,t)}
\frac{1}{m_1+m_2} \ \ g\left(t,\frac{1}{1+\frac{m_2}{m_1}}\right)\left(\frac{m_2}{m_1+m_2}\right)^{it}, \ \ t\to\infty,
\end{equation}
where the set $N$ is defined by
\begin{equation}\label{f.9}
N(\delta_1,\delta_2,t)=\left\{m_j=1,\ldots,[t], \ j=1,2,\ \frac{m_2}{m_1}\in\left( t^{1-\delta_2}-1,t^{1-\delta_1}-1\right)  \right\} .
\end{equation}
On the other hand, letting $\tau t =y$, equation \eqref{f.1} becomes
\begin{align}\label{f.10}
J(\delta_1,\delta_2,t)=\frac{1}{\sqrt{t}}& \int_{t^{\delta_1}}^{t^{\delta_2}} g\left(t,\frac{y}{t}\right) e^{itf\left(\frac{y}{t}\right)} \left|\zeta\left(\frac{1}{2}+iy\right)\right|^2 dy.
\end{align}
Estimating the above integral in the same way that $I_2$ was estimated and replacing the lhs of \eqref{f.10} by the rhs of \eqref{f.8} we find
\begin{align}\label{f.11}
\mathop{\sum\sum}_{m_1,m_2\in N(\delta_1,\delta_2,t)}
 & \frac{1}{m_1+m_2} g\left(t,\frac{1}{1+\frac{m_2}{m_1}}\right)\left(\frac{m_2}{m_1+m_2}\right)^{it}\notag\\
 &= O\left(t^{\delta_2-\frac{1}{2}}\ln t \max_{y\in[t^{\delta_1},t^{\delta_2}]}g\left(t,\frac{y}{t}\right) \right), \qquad t\to \infty,
\end{align}
where $ N(\delta_1,\delta_2,t)$ is defined in \eqref{f.9}.

In the particular case of $g(t,\tau)=1/\sqrt{\tau}$, equation \eqref{f.11} becomes 
\begin{align}\label{f.12}
\mathop{\sum\sum}_{m_1,m_2\in N(\delta_1,\delta_2,t)}
  \frac{\sqrt{m_2/m_1}}{m_2^{\frac{1}{2}-it}(m_1+m_2)^{\frac{1}{2}+it}}
 = O\left(t^{\delta_2-\frac{\delta_1}{2}} \ln t\right), \qquad t\to \infty,
\end{align}
where $ N(\delta_1,\delta_2,t)$ is defined in \eqref{f.9}.

\subsection*{Example 2}

Let
\begin{equation}\label{f.13}
f(\tau)=\tau-\tau\ln\tau.
\end{equation}
Then,
\begin{align*}
F_\tau=-\ln\tau +\ln\lambda\\
F_{\tau\tau}=-\frac{1}{\tau}.
\end{align*}
The stationary point denoted by $\xi$ satisfies
\begin{equation}\label{f.14}
\xi=\lambda,\qquad  t^{\delta_1-1} \leq \lambda \leq  t^{\delta_2-1}.
\end{equation}
Then, in analogy with \eqref{f.11} we now have
\begin{align}\label{f.15}
\mathop{\sum\sum}_{m_1,m_2\in M(\delta_1,\delta_2,t)}
 & \frac{1}{m_1} g\left(t,\frac{m_2}{m_1}\right)e^{i\frac{m_2}{m_1} t}\notag\\
 &= O\left( t^{\delta_2-\frac{1}{2}} \ln t \max_{y\in[t^{\delta_1},t^{\delta_2}]}g\left(t,\frac{y}{t}\right) \right), \qquad t\to \infty,
\end{align}
where the set $M$ is defined by
\begin{equation}\label{f.16}
M(\delta_1,\delta_2,t)=\left\{m_j=1,\ldots,[t], \ j=1,2,\ \frac{m_2}{m_1}\in\left( t^{\delta_1-1},t^{\delta_2-1}\right)  \right\} .
\end{equation}

In the case of $g(t,\tau)=\frac{1}{\sqrt{\tau+\alpha(t,\tau)}}$, where $\alpha(t,\tau)$ is a given positive function (non-strictly) increasing with respect to $\tau$, equation \eqref{f.15} becomes
\begin{align}\label{f.17gen}
\mathop{\sum\sum}_{m_1,m_2\in M(\delta_1,\delta_2,t)}
  &\frac{1}{\sqrt{m_1}\sqrt{\alpha\left(t,\frac{m_2}{m_1}\right) \hspace*{0.1mm}  m_1+ m_2}} e^{i\frac{m_2}{m_1} t}\notag\\
  &= O\left( t^{\delta_2}\frac{1}{\sqrt{t^{\delta_1} + t \ \alpha(t,t^{\delta_1-1})}} \ln t \right), \qquad t\to \infty,
\end{align}
where $M(\delta_1,\delta_2,t)$ is defined in \eqref{f.16}.

In the particular case of $\alpha= 1$, namely $g(t,\tau)=1/\sqrt{\tau+1}$, equation \eqref{f.15} becomes
\begin{align}\label{f.17}
\mathop{\sum\sum}_{m_1,m_2\in M(\delta_1,\delta_2,t)}
  \frac{1}{\sqrt{m_1}\sqrt{m_1+ m_2}} e^{i\frac{m_2}{m_1} t}= O\left( t^{\delta_2-\frac{1}{2}} \ln t \right), \qquad t\to \infty,
\end{align}
where $M(\delta_1,\delta_2,t)$ is defined in \eqref{f.16}.

In the particular case of $\alpha= 0$, namely $g(t,\tau)=1/\sqrt{\tau}$, equation \eqref{f.15} becomes
\begin{align}\label{f.18}
\mathop{\sum\sum}_{m_1,m_2\in M(\delta_1,\delta_2,t)}
  \frac{1}{\sqrt{m_1 m_2}} e^{i\frac{m_2}{m_1} t}= O\left( t^{\delta_2-\frac{\delta_1}{2}} \ln t \right), \qquad t\to \infty,
\end{align}
where $M(\delta_1,\delta_2,t)$ is defined in \eqref{f.16}.

Using the fact that
\begin{align*}
\left(\frac{m_1+m_2}{m_1}\right)^{it}=\left(1+\frac{m_2}{m_1}\right)^{it}=e^{i t \ln \left(1+\frac{m_2}{m_1}\right)}\sim e^{i t \frac{m_2}{m_1}}, \qquad \frac{m_2}{m_1}=o(1),
\end{align*}
we observe that under the change of variables $m_1+m_2=n$ and $m_2=m$, and for $\delta_1$ arbitrarily small, \eqref{f.17} yields \eqref{1.30}.

Using the fact that
\begin{align*}
\left(\frac{m_1-m_2}{m_1}\right)^{it}=\left(1-\frac{m_2}{m_1}\right)^{it}=e^{i t \ln \left(1-\frac{m_2}{m_1}\right)}\sim e^{-i t \frac{m_2}{m_1}}, \qquad \frac{m_2}{m_1}=o(1),
\end{align*}
we observe that  under the change of variables  $m_1=n_1+n_2$ and $m_2=n_1$, \eqref{f.18} yields \eqref{f.12} with $(n_1,n_2)\in N(\delta_1,\delta_2,t)$ defined in \eqref{f.9}.

\section{Conclusions}

The main results presented here are the following:

\subsubsection*{1. An exact integral equation satisfied by $|\zeta(s)|^2$ and an exact relation between certain double exponential sums}
Equation  \eqref{1.3} is a linear integral equation satisfied by $|\zeta(s)|^2$. 
This equation has its origin in a certain identity relating the Riemann and Hurwitz zeta functions derived in \cite{AsF}. The derivation of \eqref{1.3} is based on the use of the Plemelj formulae.

Equations \eqref{1.21} and \eqref{1.22} provides an exact relation between the sum $S_M$ defined in \eqref{defS3} appearing in the asymptotic analysis of $\tilde{I}_3$,  and the sum $S_R$ defined in \eqref{asym_R} appearing in the asymptotic analysis of $\tilde{I}_4$.

\subsubsection*{2. The derivation of  rigorous asymptotic results}

The rigorous asymptotic analysis of $I_1$ and $I_2$  is presented in section 4. The rigorous estimation of $I_1$ is straightforward. The rigorous estimation of $I_2$ is based on Atkinson's classical estimates.

\subsubsection*{3. The derivation of formal asymptotic results}
The asymptotics of the integral  $\tilde{I}_3$ defined in \eqref{1.10}, which denotes the integral obtained from $I_3$ by replacing $|\zeta|^2$ with its large $t$ asymptotics, can be obtained via standard asymptotic techniques for integrals. Indeed, the main contributions of $\tilde{I}_3$ arise from the associated stationary points (the relevant rigorous computation is straightforward), as well as from the end points. The contribution from the upper end point is rigorously computed in [FSF], where the analysis of the  case that the stationary point approaches the upper end point  is also presented. The rigorous computation of the analogous contribution of the lower end point can be obtained in a very similar manner. However, the investigation of the contribution of the transition zone, which is due to those values of $m_2/m_1$ where the stationary point approaches the end point, remains open.

  The asymptotics of $\tilde{I}_4$ defined in \eqref{1.16}, which denotes the integral obtained by replacing  $|\zeta|^2$ with its large $t$ asymptotics, can be obtained  via  novel asymptotic techniques. Indeed, it turns out that the relevant analysis give rise to an integral along the Hankel contour $H_1$ whose integrand involves two terms. Remarkably, the Hankel integral of the first term can be computed analytically, and thus one is left with the computation of the Hankel integral of the second term, denoted by $E_4$. By deforming the Hankel contour to pass over the relevant stationary point, and by employing Cauchy's theorem, it follows that $E_4$ yields a steepest descent contribution plus a contribution due to the associated residue. The investigation of the contribution of the transition zone, which is due to the case when the steepest descent point approaches the pole, remains open.

In order to obtain the rigorous justification of \eqref{1.30} the following tasks are required:
\begin{itemize}
\item The derivation of the analogue of the linear integral equation \eqref{1.3} with $|\zeta|^2$ replaced by $S_R$.
\item The analysis of the contribution to $\tilde{I}_3$ from the lower end point of integration (which is very similar to the analysis presented in \cite{FSF}) and the analysis of the contribution of the transition zone when the stationary point approaches the end point.
\item The investigation of the transition zone when the steepest descent point approaches the pole associated with $E_4$.
\item The proof that the limits of $\epsilon\to 0$ and $k\to\infty$ occurring in Proposition 6.2 can be interchanged.
\end{itemize}

The above tasks are carried out in \cite{F}.

In addition to the results mentioned above, the novel approach introduced here suggests several further developments, some of which were discussed in the introduction. Concrete illustrations of some of these further developments were presented in section 7, namely equations \eqref{f.12} and \eqref{f.17gen}. In addition, equation \eqref{1.5} provides the basis for obtaining a significant improvement of the best estimate regarding Lindel\"of's hypothesis, through a linear Volterra integral equation of second type for the Riemann zeta function \cite{F}. 
%Alternatively, a different choice of the part integrand in \eqref{1.5}, and of the intervals $\{L_j\}_1^4$, provide the estimates for a variety of Riemann-type double exponential sums; this leads to a family of variants of \eqref{R1}. 

\newpage

\begin{appendices}

\section{Numerical verification of (6.14).}
Let $t=10^7, \ \delta_3=\delta_4=\delta=\frac{1}{4}$. Let $\tilde{J}_4$ be defined by (6.10). We compute $\tilde{J}_4$ and the leading term of the rhs of (6.14) at the following four different values of $A$:  $\{2+3i,-2+3i,-2-3i,2-3i\}$. These points are in  the four different quadrants of the complex $z$-plane. The results are shown below:
\begin{center}
\begin{tabular}{ l l l }
$A=2+3i$ & lhs$=-0.1-i0.3$ & rhs$=-0.1-i0.3$  \\
$A=-2+3i$ &  lhs$=0.1-i0.3$ & rhs$=0.1-i0.3$  \\
$A=-2-3i$ &  lhs$=4.68\times 10^{13} - i 1.56\times 10^{14}$ & rhs$=4.68\times 10^{13}-i1.56\times 10^{14}$ \\
$A=2-3i$ &  lhs$=-0.25-i0.75$ & rhs$=-0.25-i0.75$ .
\end{tabular}
\end{center}
The relative errors are given by
\begin{center}
\begin{tabular}{ l l }
$A=2+3i$  & re=$\left|\frac{\text{rhs}-\text{lhs}}{\text{rhs}}\right|=1.42\times 10^{-8}$ \\
$A=-2+3i$  & re=$\left|\frac{\text{rhs}-\text{lhs}}{\text{rhs}}\right|=1.42\times 10^{-8}$ \\
$A=-2-3i$  & re=$\left|\frac{\text{rhs}-\text{lhs}}{\text{rhs}}\right|=1.5\times 10^{-9}$ \\
$A=2-3i$  & re=$\left|\frac{\text{rhs}-\text{lhs}}{\text{rhs}}\right|=1.26\times 10^{-8}$.
\end{tabular}
\end{center}

The first term of the rhs of (6.14) is dominant in all cases except for the third case  where $A$ is in the third quadrant. In this case, as expected, the dominant term is the second term of the rhs of (6.14) with  the relevant contribution growing like $e^{ t^{\delta }\left(-\frac{\pi}{2} -\arg A \right)}$, with $\left(-\frac{\pi}{2} -\arg A \right)\approx 0.588$.\\

\section{Numerical verification of (6.26).}

\begin{figure}
\begin{center}
%\begin{subfigure}
\includegraphics[scale=0.45]{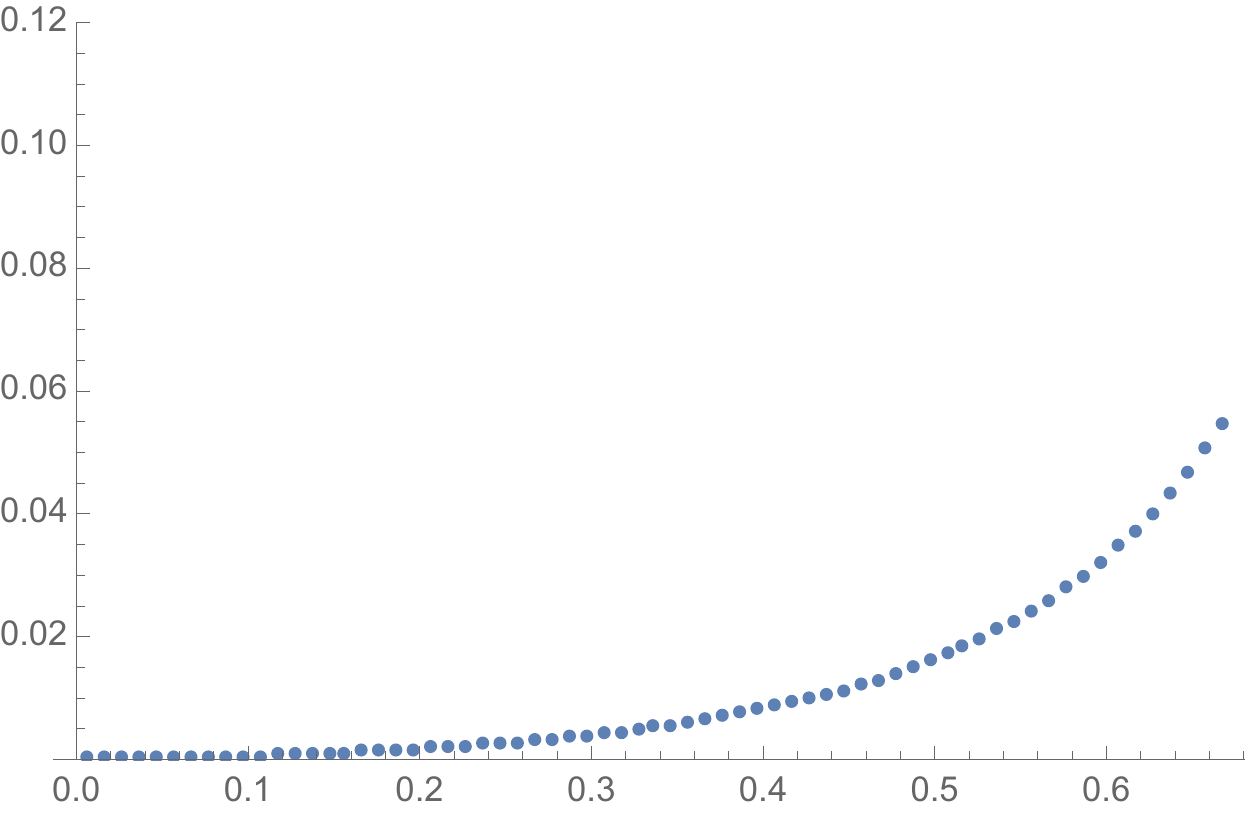}
%\end{subfigure}
%\begin{subfigure}
\qquad \includegraphics[scale=0.45]{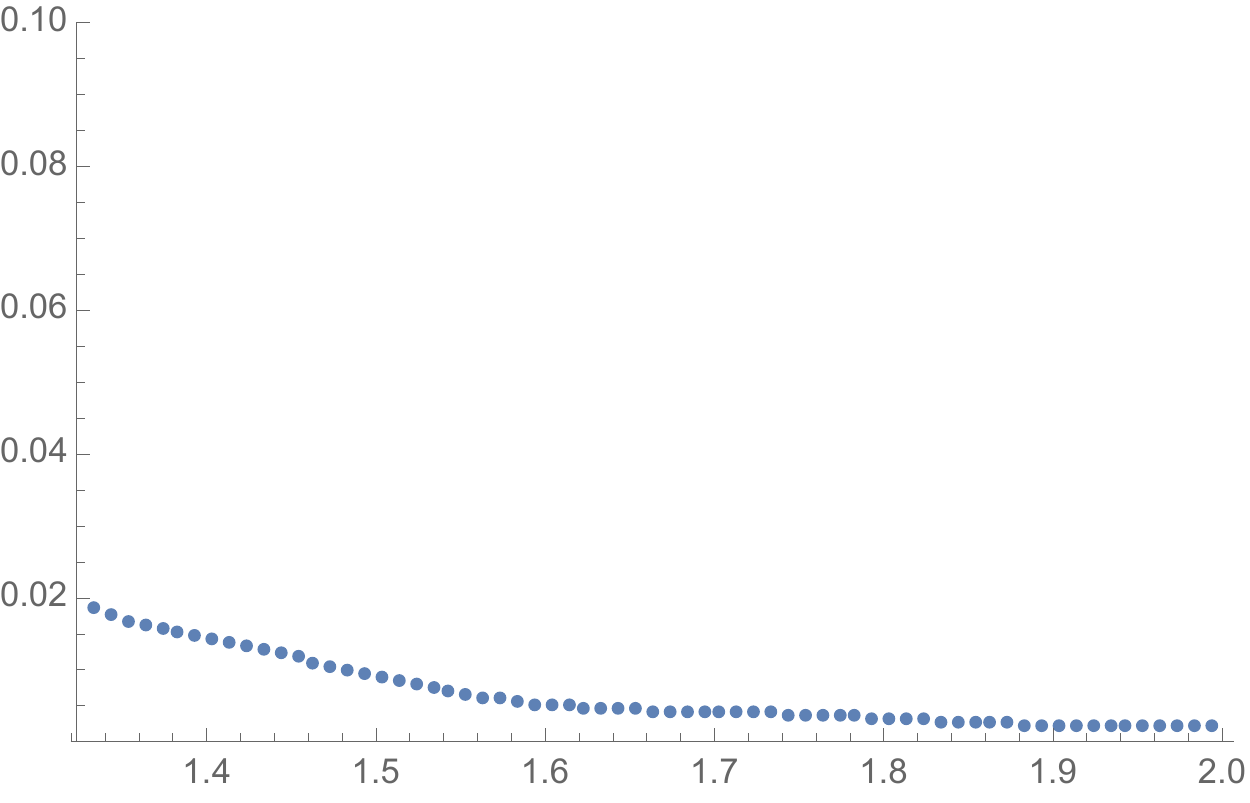}
%\end{subfigure}
\end{center}
\caption{The relative errors between the lhs and rhs of (6.26). For the left figure, the pole contribution in the rhs of (6.26) is \textit{not} taken into account, whereas for the right figure, both terms of the rhs of (6.26) are taken into consideration.}
\label{Fig-e4}
\end{figure}

Let $t=6\times10^7 + 0.45$ and $\delta_3=\frac{1}{4}$.
Let $E_4$ be defined by (6.24). Figure \ref{Fig-e4} depicts the relative error $\left|\frac{\text{rhs}-\text{lhs}}{\text{lhs}}\right|$ of (6.26) with $E_4^{SD}$ computed via (6.34). Recall the constraint (6.27), and that $M=\frac{m_1}{m_2 t}$: if $M<t^{-\delta_3}$, then the  leading asymptotic behaviour of the rhs of (6.26) is obtained by considering only the steepest descent contribution, whereas if $M>t^{-\delta_3}$ one has to consider the additional pole contribution. Figure \ref{Fig-e4} depicts the relative errors for different values of $M=a t^{-\delta_3}$: for the left figure $a\in(0,2/3)$ and the right figure $a\in(4/3,2)$. The error is small provided that the pole does not approach the stationary point, namely, $a$ does not approach the value 1.

%\section{Numerical verification of (7.20).}
%
%The mathematical analysis of $E_{34}$ is similar with $E_4$. However, it is more difficult to obtain numerical verifications for $E_{34}$ because now the interval of integration is large. Thus, we only consider a relatively small value of $t$, namely we let $t=39.17$ and $\delta_2=\frac{1}{4}$.
%Let $E_{34}$ be given by (7.19). We compute the relative error $\left|\frac{\text{rhs}-\text{lhs}}{\text{lhs}}\right|$ of (7.20) with $E_{34}^{SD}$ given by (7.25), for all points $m_j=1,\ldots,[t],\ j=1,2$, namely for $39\times 39=1521$ points. The red dots correspond to the points $(m_1,m_2)$ which satisfy the constraint (7.21), thus the computation of the rhs of (7.20) involves the contribution from pole $\omega_p$, whereas for the blue dots only the contribution from the stationary point $\omega_{sp}$ is taken into account.
%
%Again the relative error is small, unless the pole approaches the stationary point (this feature is more  striking for small values of $m_1$).
%\begin{figure}
%\begin{center}
%\includegraphics[scale=0.77]{E_34.pdf}
%\end{center}
%\caption{The relative errors between the lhs and rhs of (7.20). The blue dots correspond to the steepest descent contribution, whereas the red dots correspond to both the steepest descent and the pole contributions.}
%\label{Fig-e34}
%\end{figure}
%
%

\end{appendices}

\newpage

\section*{Erratum}
Unfortunately, there are serious numerical errors in the previous versions uploaded to arXiv, including the version that is mentioned in the announcement\footnote{https://viterbischool.usc.edu/news/2018/06/mathematician-m-d-solves-one-of-the-greatest-open-problems-in-the-history-of-mathematics/} of 2018. Fortunately, there are no errors in works that have been published. Furthermore, the approach introduced in the above versions provides a novel methodology towards attempting to prove Lindelof's hypothesis or at least towards improving dramatically the current best estimate for the large $t$-asymptotics of Riemann's zeta function.

\section*{Acknowledgement}

This project would {\it not} have been completed without the crucial contribution of Kostis Kalimeris. Kostis has studied extensively the classical techniques for the estimation of single and multiple exponential sums; these techniques are used extensively in our joint paper with Kostis \cite{KF} and some of the results of this paper are used here. Furthermore, Kostis has checked the entire manuscript and has made important contributions to the completion of some of the results presented here.

I have benefited greatly from my long collaboration with Jonatan Lenells. In particular, regarding the current work, equation \eqref{1.3} was derived in June 2015 as part of a long term collaborative project with Jonatan on the asymptotics of the Riemann zeta function and of related functions. Furthermore, the technique used  for the asymptotic analysis of $E_4$ and $E_{34}$ was introduced by Jonatan.

The rigorous estimates of the integral $J_{3}^U$ is presented in our joint paper with
Arran Fernandez and Euan Spence.

The starting point of the approach developed here is equation (2.3) which is derived in our joint paper with Anthony Ashton \cite{AsF}.

In addition to my former students Anthony, Euan and Kostis, my current
student Arran, and my former post doctoral associate Jonatan, my former students Mihalis Dimakos and Dionysis Mantzavinos have offered me generous support and assistance during the last nine years of my investigation of the asymptotics of the Riemann zeta function.

I am grateful to the late Bryce McLeod, to Eugene Shargorodsky and Bengt Fornberg,  for collaborative attempts related to the present paper, as well as to Sir Michael Atiyah, John Toland and Peter Sarnak for their encouragement.

I thank my current student  Nicholas Protonotarios for technical assistance.

Finally, I am deeply grateful to EPSRC for many years of continuous support which currently is in the form of a senior fellowship.

%\newpage

\end{document}